\documentclass[journal,twoside,web]{ieeecolor}
\usepackage{generic}
\usepackage{amsmath,amssymb,amsfonts}
\usepackage{algorithmic}
\usepackage{graphicx}
\usepackage{subfigure}
\usepackage{amsmath,bm}
\usepackage{textcomp}
\usepackage{multirow}
\usepackage{tikz}
\graphicspath{{images/}}
\usepackage{amsmath}
\usepackage{latexsym, amssymb}
\usepackage{amsbsy}
\usepackage{amsopn, amstext}
\usepackage{ifpdf,hyperref}
\usepackage{cancel, color}
\usepackage{epstopdf}
\def\BibTeX{{\rm B\kern-.05em{\sc i\kern-.025em b}\kern-.08em
		T\kern-.1667em\lower.7ex\hbox{E}\kern-.125emX}}
\def\ttimes{\,\rotatebox[]{-90}{$\ltimes$}\,}

\def\lvplus{\vec{\lplus}}

\def\J{{\bf 1}}

\usepackage[]{algorithm2e}
\usepackage{tikz}
\graphicspath{{images/}}
\usepackage{amsmath}
\usepackage{latexsym, amssymb}
\usepackage{graphicx}
\usepackage{amsmath, amsbsy}
\usepackage{amsopn, amstext}
\usepackage{ifpdf,hyperref}
\usetikzlibrary{calc,trees,positioning,arrows,chains,shapes.geometric,decorations.pathreplacing,decorations.pathmorphing,shapes,matrix,shapes.symbols}

\usepackage{tikz}
\usepackage{amsmath}
\usepackage{bm}
\usepackage{latexsym, amssymb, amsmath, amsbsy,amsopn, amstext}
\usepackage{graphicx,hyperref,booktabs}
\usepackage{cancel, color}
\usepackage{epstopdf}
\graphicspath{{images/}}

\DeclareMathOperator{\rank}{rank}
\DeclareMathOperator{\ad}{ad}
\DeclareMathOperator{\const}{const.}
\DeclareMathOperator{\Span}{Span}
\DeclareMathOperator{\Col}{Col}
\DeclareMathOperator{\Row}{Row}

\DeclareMathOperator{\argmin}{argmin}

\DeclareMathOperator{\lcm}{lcm}

\def\cal{\mathcal}

\def\ra{\rightarrow}

\def\lra{\leftrightarrow}

\def\a{\alpha}

\def\d{\delta}

\def\D{\Delta}

\def\0{{\bf 0}}

\def\J{{\bf 1}}

\def\dd{\mathrm{d}}
\def\lvplus{\vec{+}}
\newcommand{\R}{{\mathbb R}}
\newcommand{\N}{{\mathbb N}}

\newcommand{\F}{{\mathbb F}}

\def\dsum{\mathop{\sum}\limits}

\newtheorem{thm}{Theorem}[section]
\newtheorem{cor}[thm]{Corollary}
\newtheorem{dfn}[thm]{Definition}
\newtheorem{prp}[thm]{Proposition}
\newtheorem{exa}[thm]{Example}
\newtheorem{alg}[thm]{Algorithm}
\newtheorem{rem}[thm]{Remark}
\newtheorem{lem}[thm]{Lemma}

\begin{document}

\title{Observer-Based Realization of Control Systems}
\author{Daizhan Cheng,  Xiao Zhang, Zhengping Ji, Changxi Li
	\thanks{D. Cheng is with the Key Laboratory of Systems and Control, Academy of Mathematics and Systems Science, Chinese Academy of Sciences, Beijing 100190, P.R.China, e-mail: dcheng@iss.ac.cn}
    \thanks{X. Zhang is with the Department of Applied Mathematics, Hong Kong Polytechnic University, Hong Kong SAR, P.R.China, e-mail: xiaozhang@amss.ac.cn}
    \thanks{Z. Ji is with the Department of Mathematics, Friedrich-Alexander-Universit\"{a}t Erlangen-N\"{u}rnberg, 91058 Erlangen, Germany, e-mail: zhengping.ji@fau.de }
    \thanks{C. Li is with the School of Mathematics, Shandong University, Jinan 250100, P. R. China, e-mail: lichangxi@sdu.edu.cn}
	\thanks{This work is supported partly by the National Natural Science Foundation of China (NSFC) under Grant 62073315 and 62350037.}
}

\maketitle

\begin{abstract}
A novel model reduction framework for large-scale complex systems is proposed by introducing function-type dynamic control systems via the dimension-keeping semi-tensor product (DK-STP) of matrices. Utilizing bridge matrices, the DK-STP facilitates the construction of an approximate observer-based realization (OR) of a linear control system in the form of a function-type control system, where the functions serve as observers. A necessary and sufficient condition is established for the OR-system to admit exact observer dynamics. When an exact OR-system does not exist, an extended OR-system is developed by incorporating the original system’s observers into its state. Furthermore, a minimal feedback extended OR-system is constructed, and its relationship to Kalman’s minimal realization is analyzed. Finally, the proposed approach is extended to nonlinear control-affine systems.
\end{abstract}


\IEEEpeerreviewmaketitle

\section{Introduction}

Since the naissance of modern control theory, various mathematical frameworks have been proposed for modeling the dynamics of control systems, including transfer functions, phase spaces, and the behavioral approach \cite{son,wil}. Among these, the state-space representation remains the most widely adopted methodology, particularly for addressing nonlinearities.

Consider, for instance, finite-valued multi-agent systems evolving over discrete time, with Boolean networks being a prominent example. Since Kauffman introduced Boolean networks to model genetic regulatory systems \cite{kau69}, the state-space approach has found widespread applications in the analysis of Boolean (control) networks \cite{hua00} and has subsequently been extended to other finite-valued systems, such as finite games \cite{che15}.

A generic finite-valued multi-agent control system with $n$ states (nodes or agents), $m$ inputs, and $p$ outputs can be described as follows:
\begin{align}\label{1.1}
\begin{array}{l}
\begin{cases}
X_1(t+1)=f_1(X_1(t),\cdots,X_n(t),U_1(t),\cdots,U_m(t)),\\
X_2(t+1)=f_2(X_1(t),\cdots,X_n(t),U_1(t),\cdots,U_m(t)),\\
~~~~\vdots\\
X_n(t+1)=f_n(X_1(t),\cdots,X_n(t),U_1(t),\cdots,U_m(t)),\\
\end{cases}\\
~~Y_\ell(t)=h_\ell(X_1(t),\cdots,X_n(t)),\quad \ell\in [1,p],\\
\end{array}
\end{align}
where $X_i(t), U_j(t), Y_\ell(t) \in \mathcal{D}_k := \{1,2,\ldots,k\}$ for $i \in [1,n]$, $j \in [1,m]$, and $\ell \in [1,p]$ represent the system's states, inputs, and outputs, respectively. The system dynamics are given by functions $f_i: \mathcal{D}_k^{n+m} \to \mathcal{D}_k$ and $h_\ell: \mathcal{D}_k^n \to \mathcal{D}_k$.

An alternative modeling approach employs the semi-tensor product (STP) of matrices, which has proven effective in the analysis and control of finite-valued systems. This approach proceeds by identifying a discrete variable $\alpha$ with its associated vector representation $\delta_k^\alpha$ (i.e., $\alpha \sim \delta_k^\alpha$, $\alpha \in [1,k]$), and interpreting $X_i \sim x_i \in \Delta_k$. With this representation, the system in \eqref{1.1} can be equivalently reformulated as (see \cite{che11,che12}):
\begin{align}\label{1.2}
\begin{array}{l}
\begin{cases}
x(t+1)=L\ltimes u(t)\ltimes x(t),\\
y(t)=H\ltimes x(t),
\end{cases}
\end{array}
\end{align}
where $x(t):=\ltimes_{i=1}^n x_i(t)$, $u(t):=\ltimes_{j=1}^m u_j(t)$, and $y(t):=\ltimes_{\ell=1}^p y_\ell(t)$, with $\ltimes$ being the semi-tensor product.

{The equations \eqref{1.1} and \eqref{1.2} reflect fundamentally different modeling perspectives. In \eqref{1.1}, each state variable represents a local node value, whereas in \eqref{1.2}, each state encapsulates a function defined over the entire domain of node variables \cite{zx}. As an illustrative example, if the system models an evolutionary finite game, then \eqref{1.1} describes the evolution of individual strategies, while \eqref{1.2} captures the dynamics of strategy profiles or their functional representations. This conceptual difference is depicted in Fig. \ref{Fig1.1}.}

\begin{figure}
\centering
\setlength{\unitlength}{6 mm}
\begin{picture}(12,10)(-3,0)
\thicklines
{\color{red}
\put(2.5,8){\line(0,-1){8}}
\put(2.5,8){\line(1,0){1}}
\put(3.5,8){\line(0,-1){8}}
\put(4,8){\line(0,-1){8}}
\put(4,8){\line(1,0){1}}
\put(5,8){\line(0,-1){8}}
\put(5.5,8){\line(0,-1){8}}
\put(5.5,8){\line(1,0){1}}
\put(6.5,8){\line(0,-1){8}}
\put(7,8){\line(0,-1){8}}
\put(7,8){\line(1,0){1}}
\put(8,8){\line(0,-1){8}}
}
\put(-2.5,4.5){Players}
\put(-3.3,3.5){(State-type)}
\put(2,9){Profiles}
\put(4.5,9){(Function-type)}
\put(2.7,7){$z^1$}
\put(4.2,7){$z^2$}
\put(5.7,7){$z^3$}
\put(7.2,7){$\cdots$}
\put(1.2,5.7){$x_1$}
\put(1.2,4.2){$x_2$}
\put(1.2,2.7){$x_3$}
\put(1.2,1.2){$\vdots$}
\put(2.7,5.7){$a^1_1$}
\put(4.2,5.7){$a^2_1$}
\put(5.7,5.7){$a^3_1$}
\put(7.2,5.7){$\cdots$}
\put(2.7,4.2){$a^1_2$}
\put(4.2,4.2){$a^2_2$}
\put(5.7,4.2){$a^3_2$}
\put(7.2,4.2){$\cdots$}
\put(2.7,2.7){$a^1_3$}
\put(4.2,2.7){$a^2_3$}
\put(5.7,2.7){$a^3_2$}
\put(7.2,2.7){$\cdots$}
\put(2.7,1.2){$\vdots$}
\put(4.2,1.2){$\vdots$}
\put(5.7,1.2){$\vdots$}
\put(7.2,1.2){$\ddots$}
\thinlines
{\color{blue}
\put(1,6.5){\line(0,-1){1}}
\put(1,5.5){\line(1,0){8}}
\put(1,5){\line(0,-1){1}}
\put(1,4){\line(1,0){8}}
\put(1,3.5){\line(0,-1){1}}
\put(1,2.5){\line(1,0){8}}
\put(1,2){\line(0,-1){1}}
\put(1,1){\line(1,0){8}}
}
\end{picture}
\caption{State-type vs function-type dynamics\label{Fig1.1}}
\end{figure}

{Based on the nature of the variables involved in differentiation (or difference), we refer to \eqref{1.1} as a state-type system, and \eqref{1.2} as a function-type system. In the context of differential dynamic systems, a model is termed state-type if each component equation governs an individual state variable $x_i$. Conversely, a system is function-type if the equations describe the evolution of functions of the state, such as system outputs. For instance, output-based dynamic systems are function-type, since each output is a function of $x$. Function-type representations of finite-valued systems have proven to be highly effective in system analysis and control design \cite{lu17,che21,che19,li18,yan22}. In particular, when the function-type representation is designed through special choices of nodes and binary operations in finite-valued networks, it significantly reduces the model complexity \cite{li,ji}.}

We now extend the discussion to systems with continuous state spaces. Consider the classical Kalman state-space representation of a discrete-time linear control system:
\begin{align}\label{1.201}
\begin{cases}
x(t+1)=Ax(t)+Bu(t)\\
y(t)=Hx(t)
\end{cases}
\end{align}
its continuous-time counterpart:
\begin{align}\label{1.3}
\begin{cases}
\dot{x}(t)=Ax(t)+Bu(t)\\
y(t)=Hx(t)
\end{cases}
\end{align}
and a general nonlinear control-affine system:
\begin{align}\label{1.4}
\begin{cases}
\dot{x}(t)=f(x(t))+\dsum_{j=1}^mg_j(x(t))u_j(t)\\
y_{\ell}(t)=h_{\ell}(x(t))
\end{cases}
\end{align}
Naturally, we are led to ask: what are the corresponding function-type systems for these classical models? And how can they facilitate system analysis and synthesis?

This paper explores the observer-based realization (OR-system) of classical control systems as function-type representations. Here, the functions being differentiated are interpreted as observers. Thus, a function-type system derived from \eqref{1.3} or \eqref{1.4} is also referred to as an OR-system.

Several prerequisites are necessary to construct OR-systems for systems like \eqref{1.3} and \eqref{1.4}:
\begin{enumerate}
    \item Dimension-varying structure: Since the dimension of outputs may differ from that of the state, we consider such systems to be dimension-varying. The state space is embedded in the union of Euclidean spaces: $$
\R^{\infty}:=\cup_{n=1}^{\infty}\R^n.$$ To work within this space, we must define its vector space structure, topology, and cross-dimensional projection \cite{che19}. 

\item Semi-tensor product: The STP extends the conventional matrix product to matrices of arbitrary dimensions while preserving essential algebraic properties. First introduced in \cite{che01}, it has found numerous applications in Boolean networks \cite{lu17}, finite games \cite{che21}, dimension-varying systems \cite{che19}, finite automata \cite{yan22}, and coding theory \cite{zho18}. A recent development is the dimension-keeping STP (DK-STP) \cite{che24}, which equips the set of $m \times n$ matrices with a unitary ring structure and allows defining analytic functions for non-square matrices.

\item Construction of OR-systems: We begin with the discrete-time linear system \eqref{1.201}, whose corresponding state-observer (SO-) system is
\begin{align}\label{1.401}
y(t+1)=HAx(t)+HBu(t):=Mx(t)+Nu(t).
\end{align}
By replacing $M x(t)$ with $M \ttimes y(t)$, we obtain the following OR-system: $M\ttimes y(t)$ as
\begin{align}\label{1.402}
y(t+1)=M\ttimes y(t)+Nu(t),
\end{align}
where ``$\ttimes$" denotes the DK-STP.\\
Similarly, for the continuous-time case \eqref{1.3}:
\begin{align}\label{1.5}
\dot{y}(t)=HAx(t)+HBu(t):=Mx(t)+Nu(t).
\end{align}
which leads to the OR-system:
\begin{align}\label{1.6}
\dot{y}(t)=M\ttimes y(t)+Nu(t).
\end{align}
\end{enumerate}

\begin{figure}
\centering
\setlength{\unitlength}{5 mm}
\begin{picture}(17,9)
\thicklines
\put(0,0){\framebox(5,7){}}
\put(0.2,3){\begin{tiny}$\dot{x}(t)=
Ax(t)+Bu(t)$\end{tiny}}
\put(2,5){$\Sigma$}
\put(9,0){\framebox(8,1.5){\begin{tiny} $\dot{y}^s(t)=H^s\ttimes y^s(t)+B^su^s(t)$\end{tiny}}}
\put(9,3.5){\framebox(8,1.5){\begin{tiny}$\dot{y}^2(t)=H^2\ttimes y^2(t)+B^2u^2(t)$\end{tiny}}}
\put(9,5.5){\framebox(8,1.5){\begin{tiny}$\dot{y}^1(t)=H^1\ttimes y^1(t)+B^1u^1(t)$\end{tiny}}}
\put(5,0.75){\vector(1,0){4}}
\put(5,4.25){\vector(1,0){4}}
\put(5,6.25){\vector(1,0){4}}
\put(5.5,1){\begin{tiny}$y^s(t)=H^sx(t)$\end{tiny}}
\put(5.5,4.5){\begin{tiny}$y^2(t)=H^2x(t)$\end{tiny}}
\put(5.5,6.5){\begin{tiny}$y^1(t)=H^1x(t)$\end{tiny}}
\put(7,2.5){$\vdots$}
\put(0,7.5){Original System}
\put(10,7.5){OR-Systems}
\end{picture}
\caption{OR-systems of a large scale system \label{Fig1.2}}
\end{figure}

Fig. \ref{Fig1.2} illustrates how the OR-system framework enables complexity reduction for large-scale systems. If a particular system property, represented by a function involving only a subset of the state variables, is of interest, then an OR-system corresponding to that function can be formulated. This localized representation involves far fewer states, significantly reducing complexity. Multiple such OR-systems can coexist, each targeting different system properties.

This methodology has been explored for finite-valued systems in \cite{che23b}. In this paper, we aim to extend the OR-system framework to continuous-time control systems in Euclidean state spaces.

The main contributions of this paper are as follows:

\begin{itemize}
  \item We establish necessary and sufficient conditions under which an OR-system exactly replicates the dynamics of the original system. When this is not achievable, we introduce the extended OR-system, which includes the original observers as part of its state.
  \item We construct the minimal feedback OR-system, which achieves an exact realization with the smallest possible dimension, and compare it with Kalman's minimal realization.
  \item We extend the entire OR-system framework, originally developed for linear systems, to affine nonlinear control systems under mild regularity assumptions.
\end{itemize}

The rest of this paper is organized as follows: Section II provides  preliminaries  including the topological structure of the mixed-dimensional Euclidean space, and the DK-STP of matrices. Section III  first introduce the weighted DK-STP and then demonstrate how the matrix function of square matrices can be extended to the ring of non-square matrices under the DK-STP. The DK-STP based dynamic (control) systems are proposed and their solutions are obtained in Section IV. Section V uses the DK-STP with bridge matrices to construct approximate OR-systems. As an example, the natural OR-system structure of a singular system is revealed, while  necessary and sufficient conditions for a linear system to possess an exact OR-system are presented. Section VI investigates the OR-system of linear control systems, and the technique for constructing state-feedback exact OR-systems and state-feedback extended OR-systems is presented in details. Section VII extends the results obtained for linear control systems to affine nonlinear systems. Section VIII is a brief conclusion. The appendix provides detailed description of the computation of the largest $(A,B)$-invariant subspace contained in the kernel of the observers.

Before ending this section we give a list of notations.

\begin{itemize}

\item $\R^n$: the $n$-dimensional real Euclidean space.

\item $\R^n_*$: the dual space of $\R^n$.

\item ${\cal M}_{m\times n}$: the set of $m\times n$ matrices.

\item $\overline{\cal M}_{m\times n}$: the extended ring of $m\times n$ matrices by adding an identity.

\item $[a,b]$: the set of integers $\{a,a+1,\cdots,b\}$, where $a\leq b$.

\item $\lcm(n,p)$: the least common multiple of $n$ and $p$.

\item $\d_n^i$: the $i$-th column of the identity matrix $I_n$.

\item $\D_n:=\left\{\d_n^i\vert i=1,\cdots,n\right\}$.

\item $\J_k:=(\underbrace{1,\cdots,1}_k)^\mathrm{T}$.

\item ${\cal L}_{m\times n}$: the set of logical matrices, i.e.
${\cal L}_{m\times n}=\{ A\in {\cal M}_{m\times n}\;|\; \Col(A)\subset \D_m\}$.

\item $\Col(A)$: the set of columns of matrix $A$, and $\Col_i(A)$ is the $i$-th column of $A$.

\item $\Row(A)$: the set of rows of matrix $A$, and $\Row_i(A)$ is the $i$-th row of $A$.

\item $\Span(\cdot)$: the subspace or dual subspace generated by $\cdot$.

\item $\perp$: perpendicular relation between $x\in \R^n$ and $h\in \R^n_*$.

\item $\ltimes$: the semi-tensor product of matrices.

\item $\ttimes$: the dimension-keeping semi-tensor product of matrices.

\item $\overline{\cal V}_*^{(A)}$: the $A$-invariant closure of the dual subspace ${\cal V}_*\subset \R^n_*$.

\item $\overline{\cal V}_*^{(A,B)}$: the $(A,B)$-invariant closure of the dual subspace ${\cal V}_*\subset \R^n_*$.

\item $T(M)$: tangent space of manifold $M$.

\item $T*(M)$: cotangent space of manifold $M$.

\item $V(M)$: the set of smooth vector fields on $M$.

\item $V^*(M)$: the set of smooth covector fields on $M$.

\item $\D\subset T(M)$: distribution on $M$.

\item$\D_*\subset T^*(M)$: co-distribution on $M$.

\item ${\cal E}_*\subset T^*(M)$: exact co-distribution on $M$.

\item $\ad_f(g)=[f,g]$: the Lie bracket of $f,g\in V(M)$.
\end{itemize}

\section{Mathematical Preliminaries}

\subsection{Topology and Vector Space Structure on $\R^{\infty}$}
We begin by introducing the notion of dimension-free spaces. Consider the union of Euclidean spaces of all dimensions:
\begin{align*}
{\cal V}=\R^{\infty}:=\bigcup_{n=1}^{\infty}\R^n.
\end{align*}
This subsection shows how to endow $\mathcal{V}$ with the structure of a topological (pseudo-) vector space \cite{che19,che20}.

\begin{dfn}\label{d2.1.1} Let $x\in \R^p$, $y\in \R^q$, and $t=\lcm(p,q)$ be the least common multiple of $p$ and $q$. Then the semi-tensor addition of $x$ and $y$ is defined as follows:
\begin{align}\label{2.1.2}
x\lvplus y:=\left(x\otimes \J_{t/p}\right) + \left(y\otimes \J_{t/q}\right)\in \R^{t}.
\end{align}
\end{dfn}
Together with the conventional scalar multiplication, this addition turns $\R^{\infty}$ into a pseudo vector space \cite{abr78}.

To define a topology on $\R^{\infty}$, we first introduce an inner product, from which a norm and metric are derived.

\begin{dfn}\label{d2.1.2} Let $x\in \R^p$, $y\in \R^q$ and $t=\lcm(p,q)$. Define the inner product as:
\begin{align*}
\left<x,y\right>_{{\cal V}}:=\frac{1}{t}\left<x\otimes \J_{t/p},y\otimes \J_{t/q}\right>.
\end{align*}
Then the norm is given by $\Vert x\Vert_{{\cal V}}:= \sqrt{\left<x,x\right>_{{\cal V}}}$, and the distance between $x$ and $y$ is $d_{{\cal V}}(x,y):=\Vert x-y\Vert_{{\cal V}}$.
\end{dfn}

The topology on $\R^{\infty}$, denoted by ${\cal T}_d$, is generated by the metric $d_{{\cal V}}$.

\begin{rem}\label{r2.1.201}
\begin{enumerate}
\item Let $x,y\in \R^{\infty}$. $x,y$ are said to be equivalent, denoted by $x\lra y$, if there exist $\J_i$ and $\J_j$ such that $x\otimes \J_i=y\otimes \J_j$. Equivalently, $x\lra y$ if and only if $d_{{\cal V}}(x,y)=0$.
\item $\R^{\infty}$ with the classical scalar multiplication and addition defined by (\ref{2.1.2}) satisfies all vector space requirements except that the zero element is not unique. All ${\bf 0}_n=(0,\cdots,0)^\mathrm{T}\in \R^n$ that are equivalent are zero. It follows that the inverse of $x$ is
$\{y\;|\; y\lra -x\}$. Hence $\R^{\infty}$ is called a pseudo vector space.
\item Denote by
$\bar{x}:=\{y\;|\;y\lra x\}$, $r\bar{x}:=\overline{rx}$, and $\bar{x}+ \bar{y}:=\overline{x\lvplus y}$,
then the quotient space $\Omega:=\R^{\infty}/\lra$ is a vector space.
\item Under the quotient topology, $\Omega$ is a Hausdorff space. Furthermore, let $\langle\bar{x},\bar{y}\rangle_{{\cal V}}:=\langle x,y\rangle_{{\cal V}}$.
Then $\Omega$ is an inner product space, but not a Hilbert space.
\end{enumerate}
\end{rem}

Next, we consider the projection from $\R^m$ to $\R^n$.

\begin{dfn} \label{d2.1.3}  Let $\xi\in \R^m$. The projection of $\xi$ onto $\R^n$, denoted by $\pi^m_n(\xi)$, is defined as
\begin{align}\label{2.1.6}
\pi^m_n(\xi):=\underset{x\in \R_n}{\argmin}\Vert \xi - x\Vert_{{\cal V}}.
\end{align}
\end{dfn}

The projection can be expressed in matrix form as follows.

\begin{prp}[\cite{che19,che23c}] \label{p2.1.4}  Let $\xi\in \R^m$ and $x_0=\pi^m_n(\xi)\in \R^n$, $t=\lcm(m,n)$. Then
\begin{align}\label{2.1.701}
x_0=\Pi^m_n\xi \in \R^n,\quad \xi\in \R^m,
\end{align}
where
 \begin{align}\label{2.1.8}
\Pi^m_n=\frac{n}{t} \left( I_n\otimes \J^\mathrm{T}_{t/n}\right)  \left( I_m\otimes \J_{t/m}\right),
\end{align}
 Moreover, $x_0\bot\left(\xi- x_0 \right)$.
\end{prp}

For background on topological spaces and functional analysis, see \cite{cho66} and \cite{con85}.

\subsection{From the Classical STP to the DK-STP}
As mentioned in the introduction, the classical STP, proposed more than two decades ago, has found wide applications. Unfortunately, it cannot be used for linear mappings over $\R^{\infty}$, because the (classical) STP between a matrix and a vector is generally not a vector. To address this, matrix-vector and vector-vector STPs have been defined \cite{che19}. A natural question is whether it is possible to define an STP that is compatible with the conventional matrix-matrix and matrix-vector products. The newly proposed DK-STP answers this question affirmatively. The concepts and results in this subsection are mainly taken from \cite{che24}.

\begin{dfn}\label{d2.3.1} 
Let $A\in {\cal M}_{m\times n}$ and $B\in {\cal M}_{p\times q}$, $t=\lcm(n,p)$. The DK-STP of $A$ and $B$, denoted by $A\ttimes B\in {\cal M}_{m\times q}$, is defined as follows.
\begin{align}\label{2.3.1}
A\ttimes B:=\left(A\otimes \J^\mathrm{T}_{t/n}\right)\left(B\otimes \J_{t/p}\right).
\end{align}
\end{dfn}

\begin{prp}\label{p2.3.2} Let $A\in {\cal M}_{m\times n}$ and $B\in {\cal M}_{p\times q}$, $t=\lcm(n,p)$. Then $A\ttimes B:=A\Psi_{n\times p}B$,
where $\Psi_{n\times p}:=\left(I_n\otimes \J^\mathrm{T}_{t/n}\right)\left(I_p\otimes \J_{t/p}\right)\in {\cal M}_{n\times p}$ is called a bridge matrix of dimension $n\times p$.
\end{prp}

Similar to the classical STP, the DK-STP satisfies the following fundamental properties:

\begin{prp}\label{p2.3.3}
\begin{itemize}
\item[(i)] When $n=p$, $A\ttimes B=AB$.
\item[(ii)] (Associativity) $A\ttimes (B\ttimes C)=(A\ttimes B)\ttimes C.$
\item[(iii)] (Distributivity) Assume $A,B\in {\cal M}_{m\times n}$ and $C$ is of arbitrary dimension. Then
\begin{align}\nonumber
\begin{array}{l}
(A\pm B)\ttimes C=A\ttimes C\pm B\ttimes C,\\
C\ttimes (A\pm B)=C\ttimes A\pm C\ttimes B.\\
\end{array}
\end{align}

\end{itemize}
\end{prp}

The following proposition characterizes the basic properties of the DK-STP.

\begin{prp}\label{p2.3.4}

Let $A, B\in {\cal M}_{p\times q}$. Then $A\ttimes B\in {\cal M}_{p\times q}$.
$\left({\cal M}_{m\times n},+, \ttimes\right)$ is a non-commutative ring without identity.
Consider ${\cal M}_{m\times n}$ as a vector space over $\R$, then $\left({\cal M}_{m\times n}, \ttimes\right)$ is an algebra over $\R$.
Define the Lie bracket as
$$
[A,B]_{\ttimes}:=A\ttimes B-B\ttimes A,\quad A,B\in {\cal M}_{m\times n},
$$
then $\left({\cal M}_{m\times n}, \ttimes\right)$ is a Lie algebra.
\end{prp}

For concepts and fundamental properties of rings and (Lie) algebras, refer to standard textbooks such as \cite{hun74}, \cite{boo86}.

\section{Analytic Functions of Non-square Matrices}\label{S3}

\subsection{Weighted and $\Psi$-Based DK-STPs}

\begin{dfn}\label{d3.1.1} A set of vectors is called a weight, denoted by ${\bf \xi}=\{\xi_n\in \R^n\;|\; n=1,2,\cdots\}$,
if (i) $\xi_1=1$; (ii) $\xi_i=(\xi_i^1,\cdots,\xi_i^i)^\mathrm{T}\neq {\bf 0}$, $\xi_i^j\geq 0$, $i\geq 2$, $j\in [1,i]$.
\end{dfn}

\begin{dfn}\label{d3.1.2}
 Let ${\bf \xi}$ and ${\bf \eta}$ be two weights,  $A\in {\cal M}_{m\times n}$, $B\in {\cal M}_{p\times q}$, and $t=\lcm(n,p)$. Then the weighted DK-STP with respect to ${\bf \xi}$ and ${\bf \eta}$ is defined as follows.
\begin{align}\label{3.1.2}
A\ttimes^{{\bf \xi}}_{{\bf \eta}}  B:=\left(A\otimes \xi^\mathrm{T}_{t/n}\right)\left(B\otimes \eta_{t/p}\right).
\end{align}
\end{dfn}

\begin{prp}\label{p3.1.3}   Let ${\bf \xi}$ and ${\bf \eta}$ be two weights,  $A\in {\cal M}_{m\times n}$, $B\in {\cal M}_{p\times q}$, and $t=\lcm(n,p)$. Then the weighted DK-STP with respect to ${\bf \xi}$ and ${\bf \eta}$ can be calculated by
\begin{align}\label{3.1.3}
A\ttimes^{{\bf \xi}}_{{\bf \eta}} B:=A\Psi^{{\bf \xi}}_{{\bf \eta}} B,
\end{align}
where the corresponding bridge matrix
\begin{align}\label{3.1.4}
\Psi^{{\bf \xi}}_{{\bf \eta}}=\left(I_n\otimes \xi^\mathrm{T}_{t/n}\right)\left(I_p\otimes \eta_{t/p}\right)\in {\cal M}_{n\times p}.
\end{align}
\end{prp}

\begin{rem}\label {r3.1.4}
It is straightforward to verify that Propositions \ref{p2.3.3} and \ref{p2.3.4} also hold for the weighted DK-STP. The design of the bridge matrix (including the choice of weight vectors) is a key issue in applications. Here, we present two approximate approaches for the design: cross-dimensional projection and the Moore-Penrose inverse. All arguments concerning exact OR-realization focus on the design of appropriate bridge matrices for particular systems.
\end{rem}

Note that if we choose ${\bf \xi}=\{\xi_n=\frac{1}{n}\J_n\;|\;n=1,2,\cdots\}$ and  ${\bf \eta}=\{\eta_n=\J_n\;|\;n=1,2,\cdots\}$, then
\begin{align}\label{3.1.5}
\left[\Psi_p\right]_{n\times m}:=\left(\Psi^{{\bf \xi}}_{{\bf \eta}}\right)_{n\times m}=\Pi^m_n,
\end{align}
where $\Pi^m_n$ is the projection matrix, defined in (\ref{2.1.8}). This matrix is called the projecting bridge matrix.

Equation (\ref{3.1.5}) reveals the relationship between cross-dimensional projection and the DK-STP from the perspective of observer-based realization. The observer-based realization searches for the best projection from observer space $\R^m$ to the state space $\R^n$, which is realized by the bridge matrix $\Psi$. By definition, the projection from $\R^m$ to $\R^n$ finds the point in $\R^n$ that is closest to the original point in $\R^m$.

Regardless of the original geometric meaning of the DK-STP, we can formally define more DK-STPs by different choices of bridge matrices.

\begin{dfn}[\cite{che24}]\label{d3.1.5}
The set $\{\Psi_{m,n}\!\in \!{\cal M}_{m\times n}|m,n>0\}$ is called a set of bridge matrices, if
\begin{itemize}
\item[(i)] $\rank(\Psi_{m,n})=\min(m,n),\quad m,n=1,2,\cdots.$
\item[(ii)] $\Psi_{n,n}=I_n, \quad n=1,2,\cdots.$
\end{itemize}
Given a set of bridge matrices, the corresponding DK-STP, denoted by $\ttimes_{\Psi}$, can be defined as follows: Let $A\in {\cal M}_{m\times n}$ and $B\in {\cal M}_{p\times q}$. Then the corresponding DK-STP is defined by
\begin{align}\label{3.1.10}
A\ttimes_{\Psi} B:=A \Psi_{n\times p} B.
\end{align}
\end{dfn}

\vskip 2mm

\noindent{\bf Assumption 1}: For convenience, we assume hereafter that
$\Psi=\Psi_p$ and $\ttimes=\ttimes_{\Psi_p}$,
where $\Psi_p$ is defined in (\ref{3.1.5}).

\subsection{Analytic Functions over Extended Rings}

Since $\left({\cal M}_{m\times n}, +, \ttimes\right)$ is a ring without identity, to express analytic functions of  $A\in {\cal M}_{m\times n}$, we introduce an artificial identity element.

\begin{dfn}\label{d3.2.1} Define an identity element, denoted by $I_{m\times n}$, such that
\begin{align}\label{3.2.1}
A\ttimes I_{m\times n}=I_{m\times n}\ttimes A=A,\quad \forall A\in {\cal M}_{m\times n}.
\end{align}
\begin{align}\label{3.2.2}
I_{m\times n}\ttimes x:=I_m\ttimes x,\quad \forall x\in \R^{\infty}.
\end{align}
\end{dfn}

\begin{rem}\label{r3.2.101}
$I_{m\times n}$ is not a matrix, that is to say, it has no matrix expression. Equation (\ref{3.2.1}) defines the product of $I_{m\times n}$ with any element in the ring ${\cal M}_{m\times n}$, while equation (\ref{3.2.2}) defines its action on $\R^{\infty}$.
\end{rem}


\begin{rem}\label{r3.2.102}
Using (\ref{3.2.2}), we have the following equations for the action of $I_{m\times n}$:
\begin{itemize}
\item[(i)] Let $x\in \R^p$, then
\begin{align}\label{3.2.201}
\begin{array}{l}
I_{m\times n}\ttimes x=I_m\ttimes x=I_m\Psi_{m\times p}x=\Psi_{m\times p}x.
\end{array}
\end{align}
\item[(ii)] Similarly as in the case where $A\in {\cal M}_{m\times n}$ in (\ref{3.2.1}), one may define
 \begin{align}\label{3.2.202}
I_{m\times n}\ttimes x=I_{m\times n}\Psi_{n\times p}x, \quad  x\in \R^{p}.
\end{align}
Comparing it with (\ref{3.2.201}), we have $I_{m\times n}\Psi_{n\times p}=\Psi_{m\times p}$.
\item[(iii)] Assume $x\in \R^n$. Similar to (\ref{3.2.201}), we have
\begin{align}\label{3.2.204}
I_{m\times n} x=\Psi_{m\times n}x, \quad  x\in \R^{n}.
\end{align}
\end{itemize}
\end{rem}

\begin{prp}\label{p3.2.3}
$\left(\overline{\cal M}_{m\times n}, +,\ttimes\right)$ is a ring with identity $I_{m\times n}$, where

$$
\begin{array}{ccl}
\overline{\cal M}_{m\times n}&:=&{\cal M}_{m\times n}\bigcup \{rI_{m\times n}\;|\;r\in \R\}\\
~&=&\{rI_{m\times n}+A\;|\; r\in \R,\;A\in {\cal M}_{m\times n}\},
\end{array}
$$
$$
\begin{array}{l}
(r_1I_{m\times n}+A)+(r_2I_{m\times n}+B):=\\(r_1+r_2)I_{m\times n}+(A+B),
\quad  A,B\in {\cal M}_{m\times n},
\end{array}
$$
and
$$
\begin{array}{l}
(r_1I_{m\times n}+A)\ttimes (r_2I_{m\times n}+B)\\
~~:=(r_1r_2)I_{m\times n}+(r_1B+r_2A+A\ttimes B).
\end{array}
$$
\end{prp}

\begin{rem}\label{r3.2.301}
$\left(\overline{\cal M}_{m\times n}, +\right)$ is a vector space of dimension $mn+1$.
Equation~(\ref{3.2.202}) means that the action of $I_{m\times n}$ on $\R^{\infty}$ depends on the definition of $\ttimes$. More precisely, it depends on the corresponding bridge matrix $\Psi$, which determines $\ttimes$. Specifically, we may replace $I_{m\times n}$ by $I_{m\times n}\Psi$.
\end{rem}

Let $A\in {\cal M}_{m\times n}$. Define
\begin{align*}
A^{\langle k\rangle }:=
\begin{cases}
I_{m\times n},\quad k=0,\\
\underbrace{A\ttimes \cdots \ttimes A}_k,\quad k\geq 1.
\end{cases}
\end{align*}

Let $p(x)=c_nx^n+c_{n-1}x^{n-1}+\cdots+c_1x+c_0$ be a polynomial with $c_i\in \R$, $i\in [0,n]$, $A\in {\cal M}_{m\times n}$. Define the polynomial $p$ of $A$ as
\begin{align*}
\begin{array}{l}
p\langle A\rangle:=c_nA^{\langle n\rangle }+c_{n-1}A^{\langle n-1\rangle }+\cdots+c_1A+c_0I_{m\times n}\\
~~~~~\in
\overline{ {\cal M}}_{m\times n}.
\end{array}
\end{align*}

Analytic functions of non-square matrices can be defined using the above expression and the Taylor series expansions of analytic functions.

\begin{dfn}\label{d3.2.5} Given $A\in {\cal M}_{m\times n}$, one has
\begin{align*}
    &e^{\langle A\rangle }:=I_{m\times n}+A+\frac{1}{2!}A^{\langle 2\rangle }+\frac{1}{3!}A^{\langle 3\rangle }+\cdots.\\
&\cos\langle A\rangle :=I_{m\times n}-\frac{1}{2!}A^{\langle 2\rangle }+\frac{1}{4!}A^{\langle 4\rangle }+\cdots.\\
&\sin\langle A\rangle :=A-\frac{1}{3!}A^{\langle 3\rangle }+\frac{1}{5!}A^{\langle 5\rangle }+\cdots.\\
&\ln\langle I_{m\times n}+A\rangle:=A-\frac{1}{2}A^{\langle 2\rangle }+\frac{1}{3}A^{\langle 3\rangle }+\cdots.\\
&(I_{m\times n}+A)^{\langle\a\rangle}:=I+\a  A+\frac{\a (\a-1)}{2!}A^{\langle 2\rangle }+\\
&\qquad\qquad\qquad\quad~\frac{\a (\a-1) (\a-2)} {3!}A^{\langle 3\rangle }+\cdots.\\
&\cosh\langle A\rangle :=I_{m\times n}+\frac{1}{2!}A^{\langle 2\rangle }+\frac{1}{4!}A^{\langle 4\rangle }+\cdots.\\
&\sinh\langle A\rangle :=A+\frac{1}{3!}A^{\langle 3\rangle }+\frac{1}{5!}A^{\langle 5\rangle }+\cdots.
\end{align*}

\end{dfn}

All other analytic functions of a non-square matrix $A$ can be defined in a similar manner using their Taylor series expansions.

To determine when these analytic functions are well-defined, it is necessary to understand the conditions under which the corresponding Taylor series converge. For this reason, it is essential to define an appropriate norm for matrices.

\begin{dfn}\label{d3.2.6} Let $A\in {\cal M}_{m\times n}$. The DK-norm of $A$, denoted by $\Vert A\Vert_{\ttimes}$, is defined as
\begin{align*}
\Vert A\Vert_{\ttimes}:=\sup_{{\bf 0}\neq x\in \R^m}\frac{\Vert A\ttimes x\Vert_{{\cal V}}}{\Vert x\Vert_{{\cal V}}}.
\end{align*}
\end{dfn}

According to Definition \ref{d3.2.6}, it follows that
\begin{align}\label{3.2.14}
\Vert A\Vert_{\ttimes}=\sqrt{\sigma_{\max}\left(\Psi_{n\times m}^{\mathrm{T}}A^{\mathrm{T}}A\Psi_{n\times m}\right)}.
\end{align}
The convergence of the Taylor series expansions of analytic functions of non-square matrices can be verified using (\ref{3.2.14}).

\section{DK-STP Based Dynamic Systems}
\subsection{Action of Matrices on $\R^{\infty}$}

Denote the set of all matrices by
$${\cal M}:=\bigcup_{m=1}^{\infty}\bigcup_{n=1}^{\infty}{\cal M}_{m\times n}.$$

\begin{dfn}\label{d4.1.1} Let $A\in {\cal M}$ and $x\in \R^{\infty}$. The action of ${\cal M}$ on $\R^{\infty}$, denoted by $\pi_A:\R^{\infty}\ra \R^{\infty}$, is defined as
\begin{align*}
\pi_A(x):=A\ttimes x.
\end{align*}
\end{dfn}

\begin{rem}\label{r4.1.2} Recall Assumption 1,
assume $A \in {\cal M}_{m\times n}$, $x\in \R^p$, and $t=\lcm(n,p)$, then
\begin{align*}
\pi_A(x):=A\Psi_{n\times p} x=\frac{n}{t}A\left(I_n\otimes \J^\mathrm{T}_{t/n}\right)\left(I_p\otimes \J_{t/p}\right)x.
\end{align*}
\end{rem}

The mapping $\pi_A$ has the following properties.

\begin{prp}\label{p4.1.3} Let $A,B\in {\cal M}$ and $x\in \R^{\infty}$. Then $$A\ttimes(B \ttimes x)=(A\ttimes B)\ttimes x.$$
\end{prp}

\begin{prp}[\cite{che24}]\label{p4.1.4}
Let $A\in {\cal M}_{m\times n}$.
\begin{itemize}
\item[(i)] $\R^m$ is an $A$-invariant subspace. That is,
$A\ttimes x\in \R^m$, $\forall x\in \R^m$.
\item[(ii)] There exists a unique $\Pi_A\in {\cal M}_{m\times m}$, such that
$A\ttimes x=\Pi_Ax$, $\forall x\in \R^m$.
\item[(iii)]$\Pi_A=A\Psi_{n\times m}.$
\end{itemize}
\end{prp}

The construction of polynomial functions of non-square matrices implies the following generalized Cayley-Hamilton theorem:
\begin{thm}[\cite{che24}]\label{t4.1.5}
 Let $p(x)=x^m+c_{m-1}x^{m-1}+\cdots+c_1x+c_0$ be the characteristic function of $\Pi_A$, then
\begin{align*}
A^{\langle m+1\rangle}+c_{m-1}A^{\langle m\rangle}+\cdots+p_0A=0.
\end{align*}
\end{thm}

\begin{rem}\label{r4.1.6} In Theorem \ref{t4.1.5}, if $n<m$ then the characteristic function of $\Pi_{A^\mathrm{T}}$ can be used to reduce the degree of the characteristic function.
\end{rem}

\subsection{Quasi-Dynamic Systems Over ${\cal M}_{m\times n}$}

The following definition is based on those given in \cite{ahs87}, \cite{jos05}, and \cite{liu08}.

\begin{dfn}\label{d4.2.1}
\begin{itemize}
\item[(i)] Let ${\cal G}=(G,*)$ be a monoid acting on a set $X$. $({\cal G},X)$ is called a semigroup system (briefly, S-system) if there exists a mapping $\pi:G \times X\rightarrow X$, denoted by $(g,x)\mapsto gx$, satisfying
\begin{itemize}
\item[(a)] $g_1(g_2x)=(g_1*g_2)x$, $g_1,g_2\in G$, $\forall x \in X$.
\item[(b)] Let $e\in G$ be the identity, then $e(x)=x$, $\forall x\in X$.
\end{itemize}
If only condition (a) holds, $(G,X)$ is called a quasi-S-system.

\item[(ii)] Let $X$ be a topological space, and $(G,X)$ is a (quasi-) S-system system.  $(G,X)$ is called a (quasi-) dynamic system, if for each $g\in G$ the $\pi_g:X\ra X$, defined by $\pi_g: x\mapsto gx$, $x\in X$,
is continuous.
\end{itemize}
\end{dfn}

Consider the action of monoid $\left(\overline{ {\cal M}}_{m\times n},\ttimes\right)$  on $\R^{\infty}$. 

The discrete-time S-system can be constructed as
\begin{align}\label{4.2.5}
y(t+1)=A(t)\ttimes y(t),~A(t)\in  {\cal M}_{m\times n}, ~y(0)\in \R^{\infty}.
\end{align}
Similarly, the continuous-time S-system can be constructed as
\begin{align}\label{4.2.6}
\dot{y}(t)=A(t)\ttimes y(t), ~A(t)\in  {\cal M}_{m\times n}, ~y(0)\in \R^{\infty}.
\end{align}

The systems (\ref{4.2.5}) and (\ref{4.2.6}) are quasi-dynamic systems. Particularly, if the state space is restricted to $\R^m$, they become classical linear dynamic systems.

In this paper, we only focus on constant systems, that is, we assume that:

\vskip 2mm

\noindent{\bf Assumption 2}: $A(t)=A$, $\forall t\in\N$.

\vskip 2mm

Consider system (\ref{4.2.5}). A straightforward computation shows that
$$
\begin{array}{l}
y(1)=A\ttimes y_0:=y_1\in \R^m.\\
y(2)=A\ttimes y(1)=A^{\langle 2\rangle }\ttimes y_0=A^{\langle1\rangle}\Psi_{n\times m}y_1\\
y(3)=A\ttimes y(2)=A^{\langle 3\rangle }\ttimes y_0=A^{\langle 2\rangle }\Psi_{n\times m}y_1    \\
\cdots\\
\end{array}
$$
Obviously, the trajectory of (\ref{4.2.5}) is
\begin{align*}
y(t)=
    A^{\langle t\rangle } \ttimes y_0=\left(A \Psi_{n\times m}\right)^{t-1} y_1, ~t\geq 1,
\end{align*}
where $y_1=y(1)$.

\begin{rem}\label{r4.2.3}
\begin{itemize}
\item[(i)] Equation (\ref{4.2.5}) is a quasi-dynamic system, because the identity $I_{m\times n}\in \overline{\cal M}_{m\times n}$ satisfies
$$
I_{m\times n}x=I_m\ttimes x=\Psi_{m\times n}x\neq x.
$$
\item[(ii)] Assume  $y_0\in \R^m$,
Then (\ref{4.2.5}) becomes
\begin{align}\label{4.2.9}
y(t+1)=\Pi_A y(t),\quad y(0)=y_0.
\end{align}
Equation (\ref{4.2.9}) is a classical discrete-time linear system over $\R^m$.
\end{itemize}
\end{rem}

Next, we compute the trajectory of (\ref{4.2.6}). We claim that the trajectory is
\begin{align}\label{4.2.11}
y(t)=e^{\langle At\rangle }\ttimes y_0.
\end{align}

To prove the claim, we use the Taylor series expansion to express (\ref{4.2.11}) as
\begin{align*}
\begin{array}{ccl}
y(t)&=&\left(I_{m\times n}+tA+\frac{t^2}{2!}A \Psi_{n\times m}A+ \right.\\~&~& \left.\frac{t^3}{3!}A \Psi_{n\times m}A\Psi_{n\times m}A
+\cdots\right)\ttimes y_0.
\end{array}
\end{align*}
Differentiating it yields
\begin{align}\label{4.2.1201}
\begin{array}{rl}
\dot{y}(t)=&(A+t A \Psi_{n\times m}A+ \frac{t^2}{2!}  A \Psi_{n\times m}A \Psi_{n\times m}A\\~& +\frac{t^3}{3!}A \Psi_{n\times m}A\Psi_{n\times m}A\Psi_{n\times m}A
+\cdots )\ttimes y_0.
\end{array}
\end{align}
By definition of $I_{m\times n}$ we have
$$
A=A\ttimes I_{m\times n}=A\Psi_{n\times m}I_{m\times n}.
$$
Substituting this into (\ref{4.2.1201}) yields
$\dot{y}(t)=A\ttimes y(t)$,
which completes the proof of the claim.
We therefore conclude that
\begin{prp}\label{p4.2.4}
Equation (\ref{4.2.11}) is the solution of (\ref{4.2.6}).
\end{prp}

We now extend the discussion to control systems associated with the DK-STP framework. For the discrete-time case, consider
\begin{align*}
\begin{array}{l}
y(t+1)=A\ttimes y(t)+Bu(t),\\
y(0)=x_0\in \R^{n}, \; A\in  {\cal M}_{m\times n},\; B\in  {\cal M}_{m\times r}.
\end{array}
\end{align*}

A straightforward computation shows that
\begin{align}\label{4.2.14}
\begin{array}{l}
y(t+1)=\Pi_Ay(t)+Bu(t),~t\geq 1,\\
y(1)=A\ttimes y_0+Bu(0)\in \R^{m}.
\end{array}
\end{align}
When $t\geq 1$, $y(t)$ in equation (\ref{4.2.14}) is a classical discrete-time linear control system. The corresponding control problems can be solved using the classical linear control theory.

Consider the continuous-time case, we have
\begin{align*}
\begin{array}{l}
\dot{y}(t)=A\ttimes y(t)+Bu(t),\\
y(0)\in \R^{\infty}, A\in  {\cal M}_{m\times n}, B\in  {\cal M}_{m\times r}.
\end{array}
\end{align*}

A straightforward computation verifies that
\begin{align*}
y(t)=e^{\langle At\rangle }\ttimes y_0+\int_0^te^{\langle A(t-\tau)\rangle}\ttimes Bu(\tau)d\tau.
\end{align*}

Alternatively, assume $y_0\in \R^q$.  Then it is easy to verify that
\begin{align*}
y_{0_+}:=&y(0_+)=\lim_{t\ra 0_+}e^{\langle At\rangle}y_0=I_{m\times n}\ttimes y_0
=\Psi_{m\times q}y_0.
\end{align*}

Then the system can be converted to
\begin{align}\label{4.2.17}
\dot{y}(t)=\Pi_A y(t)+ Bu(t),\quad
y(0)=y_{0_+}.
\end{align}
Here, (\ref{4.2.17}) is a classical linear control system.

\begin{rem}\label{r4.2.5} When the original system is high-dimensional, exact OR-realization becomes infeasible, as it requires complete information about the entire system. Therefore, to apply OR-realization to large-scale systems, approximation is necessary. In this context, we focus on DK-STP-based approximate systems. This section provides a systematic approach to these systems, particularly regarding the computation of their trajectories. \end{rem}

\section{OR-Systems of Linear (Control) Systems}

\subsection{Approximate OR-Systems}

Recall the linear control system given by (\ref{1.201}) or (\ref{1.3}). Its corresponding SO-system, represented by (\ref{1.401}) or (\ref{1.5}), is not yet a well-posed dynamic (control) system, since the value of $y(t+1)$ does not explicitly depend on $y(t)$.

\begin{figure}
\centering
\setlength{\unitlength}{4.5 mm}
\begin{picture}(20,18)
\thicklines
\put(4, 17.1){A. Input-Output Transition System}
\put(2,12){\framebox(1.5,1){$H$}}
\put(5.5,15){\framebox(1.5,1){$A$}}
\put(9,12){\framebox(1.5,1){$H$}}
\put(12.5,15){\framebox(1.5,1){$A$}}
\put(16,12){\framebox(1.5,1){$H$}}
\put(2.75,9.5){\oval(2,1)}
\put(1.9,9.4){{\tiny$y(t\!-\!1)$}}
\put(2.75,12){\vector(0,-1){2}}
\put(0,15.5){\vector(1,0){1.75}}
\put(2.75,15.5){\oval(2,1)}
\put(1.9,15.4){{\tiny$x(t\!-\!1)$}}
\put(2.75,15){\vector(0,-1){2}}
\put(3.75,15.5){\vector(1,0){1.75}}
\put(9.75,15.5){\oval(2,1)}
\put(9.2,15.4){{\footnotesize$x(t)$}}
\put(9.75,9.5){\oval(2,1)}
\put(9.2,9.4){{\footnotesize$y(t)$}}
\put(7,15.5){\vector(1,0){1.75}}
\put(9.75,15){\vector(0,-1){2}}
\put(9.75,12){\vector(0,-1){2}}
\put(10.75,15.5){\vector(1,0){1.75}}
\put(14,15.5){\vector(1,0){1.75}}
\put(16.75,15.5){\oval(2,1)}
\put(15.9,15.4){{\tiny$x(t\!+\!1)$}}
\put(16.75,9.5){\oval(2,1)}
\put(15.9,9.4){{\tiny$y(t\!+\!1)$}}
\put(17.75,15.5){\vector(1,0){1.75}}
\put(16.75,15){\vector(0,-1){2}}
\put(16.75,12){\vector(0,-1){2}}
\put(4, 7.8){B. State-Output Transition System}
\put(2,3){\framebox(1.5,1){$H$}}
\put(5.5,3){\framebox(2,1){$HA$}}
\put(9,3){\framebox(1.5,1){$H$}}
\put(12.5,3){\framebox(2,1){$HA$}}
\put(16,3){\framebox(1.5,1){$H$}}
\put(2.75,0.5){\oval(2,1)}
\put(1.9,0.4){{\tiny$y(t\!-\!1)$}}
\put(2.75,3){\vector(0,-1){2}}
\put(2.75,6.5){\oval(2,1)}
\put(1.9,6.4){{\tiny$x(t\!-\!1)$}}
\put(2.75,6){\vector(0,-1){2}}
\put(3.5,6){\vector(1,-1){2}}
\put(9.75,6.5){\oval(2,1)}
\put(9.2,6.4){{\footnotesize$x(t)$}}
\put(9.75,0.5){\oval(2,1)}
\put(9.2,0.4){{\footnotesize$y(t)$}}
\put(9.75,6){\vector(0,-1){2}}
\put(9.75,3){\vector(0,-1){2}}
\put(7,3){\vector(1,-1){2}}
\put(10.5,6){\vector(1,-1){2}}
\put(14,3){\vector(1,-1){2}}
\put(16.75,6.5){\oval(2,1)}
\put(15.9,6.4){{\tiny$x(t\!+\!1)$}}
\put(16.75,0.5){\oval(2,1)}
\put(15.9,0.4){{\tiny$y(t\!+\!1)$}}
\put(16.75,6){\vector(0,-1){2}}
\put(16.75,3){\vector(0,-1){2}}
\qbezier(9.5,1)(8,3.5)(9.5,6)
\put(9.3,5.6){\vector(1,2){0.2}}
\qbezier(16.5,1)(15,3.5)(16.5,6)
\put(16.3,5.6){\vector(1,2){0.2}}
\put(8,3.4){$?$}
\put(15,3.4){$?$}
\put(0,3){\vector(1,-1){2}}
\put(17.5,6){\vector(1,-1){2}}
\end{picture}
\caption{S-System with I-O vs SO-System \label{Fig5.1}}
\end{figure}

The transition processes of the input-output system (\ref{1.201}) (or (\ref{1.3})) and its corresponding SO-system are illustrated in Fig. \ref{Fig5.1}. From this figure, it is evident that, in order to make the SO-system a properly defined dynamic system, a bridge (mapping) from $\R^p$ to $\R^n$. The SO-system equipped with such a bridge is referred to as the OR-system. We provide a precise definition below:

\begin{dfn}\label{d5.1.1}
The OR-system of (\ref{1.401}) is defined by
\begin{align*}
y(t+1)=M\ttimes y(t)+Nu(t),
\end{align*}
and the OR-system of (\ref{1.5}) is defined by
\begin{align*}
\dot{y}(t)=M\ttimes y(t)+Nu(t).
\end{align*}
\end{dfn}

Hereafter, we only consider continuous-time systems. All the following arguments can be extended to the discrete-time case.

In Section 4, we have discussed how to solve an OR-system, provided that $\ttimes$ is well defined. It is clear from aforementioned discussions that $\ttimes=\ttimes_{\Psi}$; that is, $\ttimes$ is completely determined by its bridge matrix, which corresponds precisely to the "bridge mapping" required in Fig.~\ref{Fig5.1}B. In what follows, we discuss different OR-systems with respect to various products defined in Section~\ref{S3}.

\begin{itemize}
\item Projection-Based $\ttimes$:
\end{itemize}

A natural way to choose the bridge is
$\Psi=\Psi_p$,
where $\Psi_p$ is defined by (\ref{3.1.5}), because the ``bridge" defined by $\Psi_p$ coincides with the projection $\pi^p_n:\R^p\ra \R^n$. This fact reveals the physical meaning of the relationship between the projection mapping and the DK-STP.

Then we have the following result.
\begin{prp}\label{p5.1.2} Consider an SO-system
\begin{align*}
\dot{y}(t)=M\dot{x}(t)+Nu(t),\quad y(0)=y_0\in \R^{\infty},
\end{align*}
where $y(t)\in \R^p$, $x(t)\in \R^n$. Using the  project bridge matrix $\Psi_p$, its OR-system is
\begin{align*}
\dot{y}(t)=M\Pi^p_n y(t)+Nu(t),\quad y(0_+)=I_p\ttimes y_0.
\end{align*}
\end{prp}

We present an example to illustrate the procedure.

\begin{exa}\label{e5.1.3}
Consider the following SO-system:
\begin{align*}
\begin{array}{l}
\dot{y}(t)=Mx(t)+Nu(t),\\
y(0)=(1,2,0,-2,-1,-1)^\mathrm{T}\in \R^6,
\end{array}
\end{align*}
where $
M=\left[\begin{smallmatrix}
-1&-2&3&2&-3\\
-3&3&-3&2&-3\\
\end{smallmatrix}\right],
N=\left[\begin{smallmatrix}
4\\0
\end{smallmatrix}\right]
$.

Note that
$
\Pi^2_5=\left[\begin{smallmatrix}
0.4&0.4&0.2&0&0\\
0&0&0.2&0.4&0.4
\end{smallmatrix}\right]^\mathrm{T}
$, we have the OR-system as
\begin{align*}
\begin{array}{ccl}
\dot{y}(t)&=&M\ttimes y(t)+Nu(t),\\
~&=&M\Psi_{5\times 2}y(t)+Nu(t),\\
~&:=&Ly(t)+Nu(t),
\end{array}
\end{align*}
where $
L=M\Pi^2_5=\left[\begin{smallmatrix}
-0.6&0.2\\
-0.6&-1\end{smallmatrix}\right]
$. Moreover, $
\Pi^6_2=\left[\begin{smallmatrix}
1&1&1&0&0&0\\
0&0&0&1&1&1
\end{smallmatrix}\right]
$. Then, $
y(0_+)=\Pi^6_2y_0=\left[\begin{smallmatrix}
3\\
-2
\end{smallmatrix}\right]
$.
\end{exa}

\begin{itemize}
\item Least-Square $\ttimes$:
\end{itemize}

In the above projection-based approach, we only need the SO-system, which can be obtained from observed data directly. Next we consider the case when more information is available. To be specific, suppose that the state-observer mapping $y(t)=Hx(t)$ such as in equation (\ref{1.401}) or (\ref{1.5}) is known. Then one can construct a ``best" linear mapping $x(t)=\Xi y(t)$ in the least-square sense as the ``bridge".
In this case, we have
\begin{align}\label{5.1.8}
x(t)=\Xi H x(t).
\end{align}
Since $x(t)$ could be any vector, (\ref{5.1.8}) leads to an algebraic equation
\begin{align}\label{5.1.9}
I_n=\Xi H.
\end{align}

The least square solution of (\ref{5.1.9}) is
\begin{align}\label{5.1.901}
\Xi= H^{+},
\end{align}
where $H^{+}$ is the  the Moore-Penrose inverse of $H$ \cite{wan04}.
The physical meaning of (\ref{5.1.901}) is clear: we use the pseudo-inverse of the mapping $H:x(t)\ra y(t)$, denoted by $\Xi$, as the inverse mapping $\Xi:y(t)\ra x(t)$.
Particularly, if $H$ is of full row rank,
then we have
\begin{align}\label{5.1.10}
\Xi=H^{+}=H^\mathrm{T}\left(HH^\mathrm{T}\right)^{-1}.
\end{align}

\begin{dfn}\label{d5.1.3}  Define
$\Psi_+=\Psi_{n\times p}:=H^{+}.
$
The corresponding DK-STP $\ttimes:=\ttimes_{\Psi_+}$ is called the pseudo-inverse based DK-STP, denoted by
$\ttimes_+$. 
\end{dfn}

With bridge matrix $\Psi_{n\times p}$, the corresponding DK-STP $\ttimes$ has ${\cal M}_{n\times p}$ as its invariant subspace. i.e., $\ttimes: {\cal M}_{n\times p}\times {\cal M}_{n\times p}\ra {\cal M}_{n\times p}$.
This leads to the following result.
\begin{prp}\label{p5.1.4} Consider a linear control system
\begin{align*}
\begin{array}{ccl}
\dot{x}(t)&=&Ax (t)+Bu(t),\quad x(t), \; x(0)=x_0\in \R^{n},\\
y(t)&=&Hx(t),\quad y(t)\in \R^p.
\end{array}
\end{align*}
Using pseudo-inverse bridge matrix (\ref{5.1.10}), its OR-system is
\begin{align*}
\dot{y}(t)=HA\Psi_{+}y(t)+HBu(t),\quad y(0_+)=I_p\ttimes x_0.
\end{align*}
\end{prp}

\subsection{Singular Systems}

Singular control systems are endowed with a natural OR structure, as is shown next. A singular system consists of a differential part and an algebraic part, which can be described as \cite{cam80}
\begin{align}\label{5.1.1201}
\begin{cases}
E\dot{x}(t)=Ax(t)+Bu(t),\\
Fx(t)=0,
\end{cases}
\end{align}
where $E\in {\cal M}_{r\times n}$, $F\in {\cal M}_{(n-r)\times n}$.

Assume $\Theta:=\begin{bmatrix}E\\F\end{bmatrix}$,
$y(t):=Ex(t)$, $z(t):=Fx(t)$.
The SO-system of (\ref{5.1.1201}) can be expressed by
\begin{align}\label{5.1.1202}
\dot{y}(t)=Ax(t)+Bu(t).
\end{align}
Then we try to construct $\Psi$ such that (\ref{5.1.1202}) can be expressed into OR-system as $\dot{y}(t)=A\ttimes_{\Psi} y(t)+Bu(t)$.

In fact, we have $
\begin{bmatrix}
y(t)\\
z(t)
\end{bmatrix}
=\begin{bmatrix}
E\\
F
\end{bmatrix}x(t)=\Theta x(t).$
Setting $z(t)={\bf 0}$, then the least square solution is
\begin{align}\label{5.1.1204}
x(t)\approx \Theta^+
\begin{bmatrix}
y(t)\\
{\bf 0}
\end{bmatrix}:=\Psi_+y(t).
\end{align}
It leads to an approximated OR-system of (\ref{5.1.1202}) as
\begin{align}\label{5.1.1205}
\dot{y}(t)\approx A\Psi_+ y(t)+Bu(t).
\end{align}
When $\Theta$ is invertible, (\ref{5.1.1205}) becomes an exact OR-system.

\begin{rem}\label{r5.1.401} Equation (\ref{5.1.1201}) is a simplified form of singular systems. In general, a singular linear control system can be expressed by
\begin{align*}
W\dot{x}(t)=Ax(t)+Bu(t), \quad x(t)\in \R^n,
\end{align*}
where $\rank(W)=r<n$. If we decompose it into  differential  and algebraic parts, then it becomes
\begin{align*}
\begin{cases}
E\dot{x}(t)=\tilde{A}x(t)+\tilde{B}u(t), \\
Fx(t)+Du(t)=0.
\end{cases}
\end{align*}
Then the above approach is applicable with a mild change of replacing (\ref{5.1.1204}) by
\begin{align*}
x(t)\approx \Theta^+\begin{bmatrix}
y(t)\\
-Du(t)
\end{bmatrix}:=\Psi_+y(t)+\tilde{D}u(t).
\end{align*}
The arguments for system (\ref{5.1.1201}) still hold for the general case.
\end{rem}

Next, we give a numerical example to depict singular control systems and their OR-systems.

\begin{exa}\label{e5.1.5}

Consider the singular system (\ref{5.1.1201}) with

\begin{align}\label{5.1.13}
E=\left[\begin{smallmatrix}
1&0&1&0\\
0&1&0&-1
\end{smallmatrix}\right];\quad
F=\left[\begin{smallmatrix}
1&0&-1&1\\
2&1&0&0
\end{smallmatrix}\right].
\end{align}
Set
$$\Theta=
\begin{bmatrix}
E\\
F
\end{bmatrix}=
\left[\begin{smallmatrix}
1&0&1&0\\
0&1&0&-1\\
1&0&-1&1\\
2&1&0&0
\end{smallmatrix}\right].
$$
Choosing a basis of $\Span(\Col(\Theta))$ as
$$
B=
\left[\begin{smallmatrix}
0&1&0\\
1&0&-1\\
0&-1&1\\
1&0&0
\end{smallmatrix}\right].
$$
It leads to $H=BC$, where
$$
C=
\left[\begin{smallmatrix}
2&1&0&0\\
1&0&1&0\\
2&0&0&1\\
\end{smallmatrix}\right].
$$

Then straightforward computation shows that
$$
\begin{array}{l}
H^+=C^\mathrm{T}(B^\mathrm{T}HC^\mathrm{T})^{-1}B^\mathrm{T}\\
\quad\quad=\left[\begin{smallmatrix}
0.225&-0.075&0.125&0.275\\
-0.2&0.4&0&0.2\\
0.525&-0.125&-0.375&-0.025\\
0.05&-0.35&0.25&-0.05
\end{smallmatrix}\right].
\end{array}
$$
The approximate OR-system (\ref{5.1.1205}) is obtained with
$$
\Psi_+=\left[\begin{smallmatrix}
0.225&-0.075\\
-0.2&0.4\\
0.525&-0.125\\
0.05&-0.35
\end{smallmatrix}\right].
$$

Consider the singular system (\ref{5.1.1201}) with $E$ the same as in (\ref{5.1.13}) and
\begin{align}\label{5.1.14}
F=\left[\begin{smallmatrix}
1&0&-1&1\\
1&1&0&0
\end{smallmatrix}\right].
\end{align}
Then the corresponding $\Theta$ is invertible. It follows that (\ref{5.1.1205})  becomes an exact OR-system with
$$
\Psi=
\left[\begin{smallmatrix}
1&-1&0&-1\\
1&-1&-1&-2
\end{smallmatrix}\right]^\mathrm{T}.
$$
\end{exa}

\subsection{Exact OR-Systems}

In subsection A, we construct the OR-systems by using bridge matrices $\Psi_p$ and $\Psi_+$ respectively. When the original system is of large scale, its OR-system, concerning some particular properties of the original system, may have much smaller size. In subsection B it was shown that the OR-system might be exact. A critical question is: does the dynamics of observers of the OR-system coincide with the output dynamics of the original system? This subsection will provide an analysis for this.

By the self-reflexivity of $\R^n$, any row vector can be viewed as a linear mapping from $\R^n$ to $\R$.
Let $H\in {\cal M}_{p\times n}$. The subspace
$$
{\cal H}_*:=\Span(\Row(H))\subset \R^n_*,
$$
is called a dual subspace.

\begin{dfn}\label{d5.2.1} Let $A\in {\cal M}_{n\times n}$. A dual subspace ${\cal H}_*\subset \R^n_*$ is called $A$-invariant, if ${\cal H}_*A\subset {\cal H}_*$.
\end{dfn}

The following fact is well known from linear algebra.

\begin{lem}\label{l5.2.2} Let ${\cal H}_*=\Span(\Row(H))$ be a dual subspace, where $H\in {\cal M}_{p\times n}$.
${\cal H}_*$ is  $A$-invariant, if and only if, there exists
$\Xi\in {\cal M}_{p\times p}$ such that
\begin{align}\label{5.2.3}
HA=\Xi H.
\end{align}
\end{lem}

By Lemma \ref{l5.2.2}, we have the following result.

\begin{prp}\label{p5.2.3} Consider the system (\ref{1.3}). There exists a bridge matrix such that its corresponding OR-system  represents the exact dynamics of the observers of the system (\ref{1.3}), if and only if, ${\cal H}_*$ is $A$-invariant.
\end{prp}

\noindent{\it Proof.} (Sufficiency) Suppose that (\ref{5.2.3}) holds. Then the SO-system of (\ref{1.3}) becomes
\begin{align}\label{5.2.4}
\begin{array}{ccl}
\dot{y}(t)&=&HAx(t)+HBu(t)\\
~&=&\Xi y(t)+HBu(t).
\end{array}
\end{align}
Here (\ref{5.2.4}) is an OR-system. That is, the observers of (\ref{1.3}) are the trajectories of this OR-system.

(Necessity) Suppose that there exists a bridge matrix $\Psi$ such that the DK-STP $\ttimes=\ttimes_{\Psi}$ satisfies
$$
\begin{array}{ccl}
\dot{y}(t)&=&HAx(t)+HBu(t)\\
~&=&HA\ttimes y(t)+HBu(t)\\
~&=&HA\Psi Hx(t)+HBu(t).
\end{array}
$$
Then we must have $HA\Psi H=HA$.
Setting $\Xi:=HA\Psi$ implies (\ref{5.2.3}).
\hfill $\Box$

\begin{cor}\label{c5.2.4} Suppose ${\cal H}_*$ is $A$-invariant. Then by letting
$\Xi=HAH^{+}$,
$\Psi=H^{+}$,
the corresponding OR-system is exact. Particularly, when $H$ is of full row rank, we have
\begin{align}\label{5.2.5}
\begin{array}{l}
\Xi=HAH^\mathrm{T}(HH^\mathrm{T})^{-1},\\
\Psi=H^\mathrm{T}(HH^\mathrm{T})^{-1}.
\end{array}
\end{align}
\end{cor}

\noindent{\it Proof.} Without loss of generality, we prove it by assuming $H$ is of full row rank. In this case,
$$
H^{+}=H^\mathrm{T}(HH^\mathrm{T})^{-1}.
$$

It is easy to verify that ${\cal H}_*$ is $A$-invariant if and only if $\Xi=HAH^\mathrm{T}(HH^\mathrm{T})^{-1}$ satisfies
(\ref{5.2.3}). From the proof of Proposition \ref{p5.2.3} one sees that $\Psi$ satisfies
\begin{align}\label{5.2.6}
HA\Psi=\Xi=HAH^\mathrm{T}(HH^\mathrm{T})^{-1}.
\end{align}
An obvious solution of (\ref{5.2.6}) is $\Psi=H^\mathrm{T}(HH^\mathrm{T})^{-1}$.
\hfill $\Box$

\begin{rem}\label{r5.2.5} Assume ${\cal H}_*$ is $A$-invariant, then the only choice of $\Xi$ is the one in (\ref{5.2.4}), (or (\ref{5.2.5}) when $H$ is of full row rank).
But the solution for $\Psi$ may not be unique. In fact, any matrix in the form of
$\Psi=H^{+}+G$,
with $G\in {\cal M}_{n\times p}$ in the kernel of $HA$ is also a solution, while the non-uniqueness of the solution does not affect the expression of the corresponding OR-system.
\end{rem}

\begin{exa}\label{e5.2.5}
Consider a control system
\begin{align*}
\begin{array}{l}
\dot{x}(t)=Ax(t)+Bu(t),\\
y(t)=Cx(t),
\end{array}
\end{align*}
where $
A=\left[\begin{smallmatrix}
0&-2&1&-6&-9\\
-1&-3&4&-11&-13\\
4&1&-1&10&12\\
2&1&-2&7&7\\
-1&0&1&-2&0\\
\end{smallmatrix}\right]$, $
B=\left[\begin{smallmatrix}
2\\
0\\
1\\
1\\
-1\\
\end{smallmatrix}\right]$,
$$
C=\left[\begin{smallmatrix}
1&-1&1&-2&0\\
-1&0&0&-1&-2\\
\end{smallmatrix}\right].
$$
Set
$$
\Xi:=CAC^\mathrm{T}(CC^\mathrm{T})^{-1}=
\left[\begin{smallmatrix}
0&-1\\
-1&-1
\end{smallmatrix}\right].
$$
It can be verified that $
CA=\Xi C.$
Hence ${\cal C}_*=\Span(\Row(C))$ is $A$-invariant. Then we have exact OR-system as
$$
\begin{array}{ccl}
\dot{y}(t)&=&\Xi y(t)+CBu(t)\\
&=&\left[\begin{smallmatrix}
0&-1\\
-1&-1
\end{smallmatrix}\right]y(t)+\left[\begin{smallmatrix}
1\\
-1
\end{smallmatrix}\right]u(t).
\end{array}
$$

\end{exa}

\section{Extended OR-Systems}
This section addresses the problem of constructing an exact OR-system when the output dual subspace ${\cal H}_*$ is not $A$-invariant. The basic idea is to enlarge the output dual subspace so that it satisfies the invariance requirement.

\subsection{$A$-invariant Extended OR-Systems}

\begin{dfn}\label{d6.1.1} Let ${\cal H}_*\subset \R^n_*$. The smallest dual subspace that is $A$-invariant and contains ${\cal H}_*$, is called the $A$-invariant closure of ${\cal H}_*$, denoted by $\overline{\cal H}_*^{(A)}$.
\end{dfn}

\begin{rem}\label{r6.1.2}
Let ${\cal H}^i_*$, $i=1,2$ be $A$-invariant dual subspaces. Then ${\cal H}_*^1\cap {\cal H}_*^2$ is also $A$-invariant. So $\overline{\cal H}_*^{(A)}$ exists, because
$$
\overline{\cal H}_*^{(A)}=\bigcap\{{\cal S}\;|\; {\cal S}\supset {\cal H}_*~\mbox{is $A$-invariant}\}.
$$
$\overline{\cal H}_*^{(A)}$ can be calculated by an iteration as follows:
\begin{itemize}
\item[(1)]
$
{\cal V}_*^0:={\cal H}_*.
$
\item[(2)]
$
{\cal V}_*^{k+1}={\cal V}_*^k+{\cal V}_*^kA,\quad k\geq 0.
$
\item[(3)]
When the sequence reaches $k^*$ such that
$
{\cal V}_*^{k^*+1}={\cal V}_*^{k^*}$, we have
$
\overline{\cal H}_*^{(A)}={\cal V}_*^{k^*}.
$
\end{itemize}
\end{rem}

Since $\overline{\cal H}_*^{(A)}$ is $A$-invariant, we can use it to construct an OR-system.

\begin{dfn}\label{d6.1.3}
 Consider system (\ref{1.3}). Replace its observers by $\overline{\cal H}_*^{(A)}$. The resulting OR-system is called its extended OR-system.
\end{dfn}

We present a numerical example to illustrate the construction of an extended OR-system.

\begin{exa}\label{e6.1.4}
Consider a system $\Sigma$ characterized by (\ref{1.3}) with
\begin{align}\label{6.1.2}
\begin{array}{l}
A=\left[\begin{smallmatrix}
0&1&0&0&0&0\\
0&0&1&0&0&0\\
1&0&-1&1&-2&1\\
0&0&0&0&0&0\\
0&0&0&0&1&0\\
0&0&0&0&0&0\\
\end{smallmatrix}\right],\quad
B=\left[\begin{smallmatrix}
0\\
0\\
1\\
0\\
0\\
0\\
\end{smallmatrix}\right],\\
C=\left[\begin{smallmatrix}1&0&0&0&1&0\end{smallmatrix}\right].
\end{array}
\end{align}

Let ${\cal C}_*=\Span(\Row(C))$. Using the algorithm in Remark \ref{r6.1.2}, we have
$$
\overline{\cal C}_*^{(A)}=\Span(\Row(H)),
$$
where
$
H=\left[\begin{smallmatrix}
1&0&0&0&1&0\\
0&1&0&0&1&0\\
0&0&1&0&1&0\\
1&0&-1&1&-1&1\\
-1&1&1&-1&1&-1\\
\end{smallmatrix}\right].
$

Define
$z=(z_1,z_2,z_3,z_4,z_5)^\mathrm{T}:=
Hx.
$
Then we have the OR-system with respect to $\overline{\cal C}_*^{(A)}$ as
\begin{align*}
\begin{array}{ccl}
\dot{z}(t)&=&HAx(t)+HBu(t)\\
~&:=&\tilde{A}z(t)+\tilde{B}u(t),\\
y(t)&=&z_1(t).
\end{array}
\end{align*}
where
$
\tilde{A}=
\left[\begin{smallmatrix}
0&1&0&0&0\\
0&0&1&0&0\\
0&0&0&1&0\\
0&0&0&0&1\\
0&-1&1&1&0\\
\end{smallmatrix}\right];
\tilde{B}=\left[\begin{smallmatrix}
0\\0\\1\\-1\\1\end{smallmatrix}\right].
$
\end{exa}

\vskip 2mm

\begin{rem}\label{r6.1.5}

The main difference between Kalman realization and observer-based realization is that Kalman realization preserves the input-output relation of the original system, whereas observer-based realization preserves the observers (outputs) themselves.
As mentioned before, the observer (as functions of states) is to describe certain properties of the original systems concerned. Hence for our purpose, we need to keep the observers unchanged, which shows the advantage of the OR approach over Kalman realization.

\end{rem} 

\subsection{Feedback Extended OR-Systems}

\begin{dfn}\label{d6.2.1} Consider system (\ref{1.3}). A dual space ${\cal H}_*\subset \R^n_*$ is called $(A,B)$-invariant if there exists a state feedback control $u(t)=Fx(t)$ such that ${\cal H}_*$ is invariant under $A+BF$.
\end{dfn}

The following results are immediate consequences of the definition.

\begin{prp}\label{p6.2.2} ${\cal H}_*$ is an $(A,B)$-invariant dual subspace, if and only if, there exist $F$ and $\Xi$, such that
$H(A+BF)=\Xi H$.
\end{prp}

\begin{prp}\label{p6.2.3} Consider the system (\ref{1.3}). Using bridge matrix $\Psi$ and feedback control
$u(t)=Fx(t)+v(t)$,
we have the approximated OR-system
\begin{align*}
\begin{array}{ccl}
\dot{y}(t)&=&H(A+BF)x(t) +v(t)\\
~&\approx&H(A+BF) \Psi y(t)+v(t).
\end{array}
\end{align*}

If ${\cal H}_*$ is an $(A,B)$-invariant dual subspace, then the approximated OR-system becomes an exact system $\dot{y}(t)=\Xi  y(t)+v(t)$.
\end{prp}

Next, we consider the problem of how to verify if a dual subspace ${\cal H}_*$ is $(A,B)$-invariant.

\begin{prp}\label{p6.2.4} ${\cal H}_*\subset \R^n_*$ is $(A,B)$-invariant dual subspace, if and only if,
$$
{\cal H}_*^{\perp}:=\left\{x\in \R^n\;|\; h(x)=0,\forall h\in {\cal H}_*\right\}
$$
is $(A,B)$-invariant.
\end{prp}

\noindent{\it Proof.} Assume ${\cal H}_*$ is $(A,B)$-invariant, then
$$
\begin{array}{l}
~~{\cal H}_*(A+BF){\cal H}_*^{\perp}\\
=\Span(\Row(H))(A+BF){\cal H}_*^{\perp}\\
= \Span(\Row(\Xi H)){\cal H}_*^{\perp}\subset  \Span(\Row(H)){\cal H}_*^{\perp}=0,
\end{array}
$$
that is,
$$
(A+BF){\cal H}_*^{\perp}\subset {\cal H}_*^{\perp},
$$
which shows ${\cal H}_*^{\perp}$ is $(A,B)$-invariant.

Conversely, assume ${\cal H}_*^{\perp}$ is $(A,B)$-invariant,
then $(A+BF){\cal H}_*^{\perp}\subset {\cal H}_*^{\perp}$. Hence, we also have
$$
{\cal H}_*(A+BF){\cal H}_*^{\perp}=0,
$$
that is,
$$
{\cal H}_*(A+BF)\subset \left\{{\cal H}_*^{\perp}\right\}^{\perp}={\cal H}_*
$$
which implies that ${\cal H}_*$ is $(A,B)$-invariant.
\hfill $\Box$

\begin{rem}\label{r6.2.5} It is well known that a subspace ${\cal V}\subset \R^n$ is an $(A,B)$-invariant subspace if and only if $A{\cal V}={\cal V}+{\cal B}$,
where ${\cal B}=\Span(\Col(B))$.
Hence, it serves as a straightforward verification for $(A,B)$-invariant subspace. Then Proposition \ref{p6.2.4} provides a method to check whether a dual subspace is $(A,B)$-invariant.
\end{rem}

Next, consider the general case, where ${\cal H}_*$ is not $(A,B)$-invariant. To get an exact OR-system, we try to find a smallest $\Sigma_*\subset \R^n_*$ such that ${\cal H}_*\subset \Sigma_*$ and $\Sigma_*$ is $(A,B)$-invariant.

\begin{dfn}\label{d6.2.6} The smallest $(A,B)$-invariant subspace of $\R^n_*$ containing ${\cal H}_*$ is called the $(A,B)$-invariant closure of ${\cal H}_*$, denoted by $\overline{\cal H}_*^{(A,B)}$.
\end{dfn}

 We will show the existence and uniqueness of  $(A,B)$-invariant closure of ${\cal H}_*$ by constructing it explicitly.

Consider system (\ref{1.3}), there has been a standard algorithm  to calculate the largest $(A,B)$-invariant subspace contained in ${\cal H}_*^{\perp}$ \cite{won79}, by which we have the following result.

\begin{prp}\label{p6.2.7} Consider system (\ref{1.3}). Let $W\subset \R^n$ be the  largest $(A,B)$-invariant subspace contained in ${\cal H}_*^{\perp}$. Then
\begin{align*}
\overline{\cal H}_*^{(A,B)}=W^{\perp}.
\end{align*}
 \end{prp}

 \noindent{\it Proof.} First, by $W\subset {\cal H}_*^{\perp}$, we have $W^{\perp}\supset {\cal H}_*$.
Since $W$ is an $(A,B)$-invariant subspace,  $W^{\perp}$ is an $(A,B)$-invariant dual subspace.

Second, if ${\cal V}_*\supset {\cal H}_*$ is the smallest $(A,B)$-invariant dual subspace, then ${\cal V}_*^{\perp}$ is  the largest $(A,B)$-invariant subspace contained in ${\cal H}_*^{\perp}$. Hence
${\cal V}_*^{\perp}=W$. That is,
$
{\cal V}_*=W^{\perp}.
$
\hfill $\Box$

\begin{dfn}\label{d6.2.8} Consider system (\ref{1.3}). The OR-system constructed by using $\overline{\cal H}_*^{(A,B)}$ as the observer subspace is called the feedback OR-system of (\ref{1.3}).
\end{dfn}

\begin{prp}\label{p6.2.9} The feedback OR-system of (\ref{1.3}) is an exact system, which contains full information of the dynamics of the original observers.
\end{prp}

\noindent{\it Proof.} Since $\overline{\cal H}_*^{(A,B)}$ is $(A,B)$-invariant, the feedback OR-system is exact. Since the original observers are part of the state variables of the   feedback OR-system, their dynamic process is completely recovered by the feedback OR-system.

\hfill $\Box$

Since the largest $(A,B)$-invariant subspace contained in the  kernel of $H$ (i.e., ${\cal H}_*^{\perp}$) is unique \cite{won79}, the  feedback OR-system of (\ref{1.3}) on $\overline{\cal H}_*^{(A,B)}$ is also unique. It can be regarded as a minimal exact realization of (\ref{1.3}). We compare it with the classical minimal realization \cite{kal79} of (\ref{1.3}), which is recalled below.

\begin{prp}\label{p6.2.10} Consider the system (\ref{1.3}). There exists a coordinate transformation $z=z(x)$, such that (\ref{1.3}) can be transformed into the following form, (which is then called the Kalman canonical form):
\begin{align*}
\dot{z}(t)=\tilde{A}z(t)+\tilde{B}u(t),\quad
y(t)=\tilde{C}z(t),
\end{align*}
where $z(t)=(z^1(t),z^2(t),z^3(t),z^4(t))^\mathrm{T}$,
$$
\tilde{A}=\!\left[\begin{smallmatrix}
A_{11}&0&A_{13}&0\\
A_{21}&A_{22}&A_{23}&A_{24}\\
0&0&A_{33}&0\\
0&0&A_{43}&A_{44}
\end{smallmatrix}\right];\tilde{B}=\!\left[\begin{smallmatrix}
B_{1}\\
B_{2}\\
0\\
0
\end{smallmatrix}\right];\tilde{C}=\left[\begin{smallmatrix}
   C_1&0&C_3&0
\end{smallmatrix}\right],
$$
with
$(A_{11}, B_1)$ controllable and $(C_1, A_{11})$ observable. Moreover,
\begin{align}\label{6.2.9}
\begin{array}{l}
\dot{z}^1(t)=A_{11}z^1(t)+B_1u(t),\\
y(t)=C_1 z^1(t),
\end{array}
\end{align}
is a minimal realization.
\end{prp}

To compare the minimal realization (\ref{6.2.9}) of (\ref{1.3}) with its $\overline{\cal H}_*^{(A,B)}$-based
feedback OR-system, we consider a particular system.

\begin{exa}\label{e6.2.11}
Recall Example \ref{e6.1.4}.

The largest $(A,B)$-invariant subspace contained in the kernel of ${\cal C}_*$ is (one may refer to the  Appendix  for detailed calculation)
\begin{align*}
{\cal V}=\Span\left\{ \left[\begin{smallmatrix} 0\\0\\0\\1\\0\\0\end{smallmatrix}\right],
 \left[\begin{smallmatrix} 0\\0\\0\\0\\0\\1\end{smallmatrix}\right],
 \left[\begin{smallmatrix} 1\\1\\1\\0\\-1\\0\end{smallmatrix}\right]
\right\},
\end{align*}
and under the state feedback $u(t)=Fx(t)$,
with
$F=(-1,0,2,-1,2,-1)$,
we have $(A+BF){\cal V}\subset {\cal V}$.
Then the $(A,B)$-invariant closure of ${\cal C}_*$ is
$\overline{\cal C}_*^{(A,B)}=\Span\left(\Row(V^{\perp})\right),
$
where
$$
V^{\perp}=\left[\begin{smallmatrix}
1&0&0&0&1&0\\
1&-1&0&0&0&0\\
1&0&-1&0&0&0\\
\end{smallmatrix}\right].
$$

Finally, the
feedback OR-system of (\ref{6.1.2}) with respect to  $\overline{\cal C}_*^{(A,B)}$ can be obtained as
\begin{align*}
\begin{array}{ccl}
\dot{w}(t)&=&V^{\perp}\dot{x}(t)+V^{\perp}Bv(t)\\
~&=~&V^{\perp}Ax(t)+V^{\perp}Bv\\
~&=&\left[\begin{smallmatrix}
0&1&0&0&1&0\\
0&1&-1&0&0&0\\
0&1&-1&0&0&0\\
\end{smallmatrix}\right]x(t)+\left[\begin{smallmatrix}0\\-1\\0\end{smallmatrix}\right]v(t)\\
~&=&\left[\begin{smallmatrix}
0&1&0\\
0&-1&1\\
0&-1&1\\
\end{smallmatrix}\right]w(t)+\left[\begin{smallmatrix}0\\-1\\0\end{smallmatrix}\right]v(t),\\
y(t)&=&w_1(t).
\end{array}
\end{align*}

\end{exa}

\begin{rem}\label{r6.2.12}
\begin{itemize}
\item[(i)] One can see that (\ref{6.1.2}) is of the Kalman canonical form. Hence its minimal realization is
\begin{align*}
\begin{array}{ccl}
\dot{z}(t)&=&
\left[\begin{smallmatrix}
0&1&0\\
0&0&1\\
1&0&-1
\end{smallmatrix}\right]z(t)+
\left[\begin{smallmatrix}
0\\
0\\
1
\end{smallmatrix}\right]u(t),\\
y(t)&=&\left[\begin{smallmatrix} 1&0&0\end{smallmatrix}\right]z(t).
\end{array}
\end{align*}
Although this realization keeps the input-output relation of the original system unchanged, it does not keep the output functions unchanged. So this minimal realization can not be adopted for realizing the dynamics of certain chosen functions as the exact SO-systems do.

\item[(ii)] Comparing Example \ref{e6.2.11} with
 Example \ref{e6.1.4}, one sees that the feedback  OR-system has lower dimensions than the state extended OR-system. This is in general true. Because the feedback ones use largest $(A,B)$-invariant subspace to construct its perpendicular dual subspace, which is the largest feedback invariant dual subspace. The $A$-invariant dual subspace can be considered as a special $(A,B)$-invariant dual subspace, which uses feedback control $u(t)=Fx(t)$ with $F=0$. Hence, it is in general smaller than the largest one.

\item[(iii)] The feedback OR-system may simultaneously solve other control problems. For example,  the disturbance decoupling problem (DDP). Consider a linear control system with disturbances
\begin{align*}
\dot{x}(t)=Ax(t)+Bu(t)+\dsum_{k=1}^s\xi(t), \quad y(t)=Cx(t).
\end{align*}
where $\xi_k(t)$, $k\in [1,s]$ are disturbances. The DDP is solvable, if and only if, $\xi_k\in {\cal V}$, where ${\cal V}$ is the largest $(A,B)$-invariant subspace contained in kernel of $Cx$ \cite{won79}.

Now if the DDP is solvable, then we claim that all the disturbances will be eliminated from the feedback extended OR-system. This is because
$C\xi_k=0$,
and then
$$
C(A+BF)\xi_k\subset C(A+BF){\cal V}\subset  C{\cal V}=0.
$$
Continuing this process, one sees that
$\overline{{\cal C}}_*^{(A,B)}\subset {\cal V}^{\perp}$.
Hence all the disturbances are within  ${\cal V}_*^{\perp}$. The claim follows.
Similarly, if the feedback solves the stabilization problem, i.e., $A+BF$ is a Hurwitz matrix, then the    feedback extended exact OR-system is also asymptotically stable. In fact, the feedback extended OR-system is the smallest subsystem, involving observers.
\end{itemize}
\end{rem}

\section{OR-Systems of Nonlinear Control Systems}

The arguments presented in the previous sections for OR-systems can be extended to nonlinear control systems. Consider system (\ref{1.4}) with $x(t)$ evolving on an $n$-dimensional manifold $M$.  Then (\ref{1.4}) can be considered as an expression over a coordinate chart $W$ with coordinates $x(t)\in W\subset \R^n$ and $0\in W$. Hereafter, we assume $x(t)\in W\subset M$. In case that $W=\R^n$, $x(t)$ becomes a global coordinate frame. Denote by $V(W)$ and $V^*(W)$ the sets of smooth vector fields and co-vector fields on $W$ respectively (for the preliminaries on differentiable manifolds one may refer to standard textbooks such as \cite{boo86}). 

For statement ease, we assume $0\in W$ is a regular point, that is, all the distributions and co-distributions involved below are of constant dimension locally. In \eqref{1.4}, suppose that $f(x),g_i(x)\in V(W)$ are smooth vector fields. $h_j(x)$ are smooth functions. We first recall some well-known notions from nonlinear control systems.

\begin{dfn}[\cite{boo86}] \label{d7.1}
\begin{itemize}
\item[(i)] Let $\D\subset V(W)$ be a distribution. $\D$ is called nonsingular if $\dim(\D(x))=\const$, $\forall x\in W$. $\D$ is called involutive, if for any two vector fields $f(x),g(x)\in V(W)$, $[f(x),g(x)]\in \D$.
\item[(ii)] Consider a co-vector field $\xi(x)\in V^*(W)$. $\xi(x)$ is called exact, if there exists a smooth function $h(x)$ on $W$ such that $\xi(x)=\dd h(x)$.
Denote the set of co-distributions by $\D_*\subset V^*(W)$.
\end{itemize}
\end{dfn}

\begin{dfn}[\cite{isi95}]\label{d7.2}
Let $\D(x)\subset V(W)$ be a distribution. The smallest involutive distribution containing $\D(x)$ is called the involutive closure of $\D(x)$, denoted by $\overline{\D}(x)$.
A distribution $\D(x)\subset V(W)$ is said to be $(f,g)$-invariant, if
\begin{align*}
\ad_f(\D(x))\subset \D(x)+{\cal G},
\end{align*}
where ${\cal G}=\Span(g_i(x)\;|\; i\in [1,m])$.
\end{dfn}

Recall system (\ref{1.4}). We make the following assumption.

\vskip 2mm

\noindent{\bf Assumption A3}:
$h_j(0)=0,\quad j\in [1,p]$.

\begin{rem}\label{r7.201} Assumption A3 is not essential. This is to ensure the one-to-one correspondence between a function $h$ and its differential $\dd h$, as parallel to the linear case.
\end{rem}

The SO-system of system (\ref{1.4}) can be solved as
\begin{align*}
\dot{y}(t)=L_f h(x(t))+\dsum_{i=1}^mL_{g_i}h(x(t))u_i(t).
\end{align*}

Motivated by the linear case, one can see that if there exists a smooth mapping
$\psi:\R^m\ra W$,
then an approximate OR-system can be obtained as
\begin{align}\label{7.5}
\dot{y}=L_fh(\psi(y(t)))+\dsum_{i=1}^mL_{g_i}h(\psi(y(t)))u_i(t).
\end{align}

Naturally, the mappings used for the linear case are also the candidates of $\psi$, to be specific, using the projection, we have $x(t)=\Pi^p_n y(t)$; using pseudo-inverse, we have $x(t)=[\dd h(x)]^+$.

Next, we consider the case when there exists $\psi$ such that the corresponding OR-system (\ref{7.5}) is exact.

\begin{dfn}\label{d7.3}
\begin{itemize}
\item[(i)]
A co-distribution ${\cal H}_*\subset V^*(W)$ is called an exact co-distribution, if there exist $h_j(x)$, $j\in [1,p]$, such that ${\cal H}_*=\Span\{\dd h_j(x)\:|\: j\in [1,p]\}$.
Denote the set of exact co-distributions by ${\cal E}_*\subset \D_*$.
\item[(ii)] An exact co-distribution  ${\cal H}_*=\dd H=\Span\{\dd h_j(x)\:|\: j\in [1,p]\}\subset V^*(W)$ is said to be $f$-invariant, if there exists a smooth matrix $\Xi\in {\cal M}_{p\times p}$, such that $L_f(\dd H)=\Xi \dd H$.
\end{itemize}
\end{dfn}
\begin{rem}\label{r7.301}
A co-vector field is a dual vector field in the dual space $V^*(W)$. Only when it is exact, it can be generated by the differential of a function $h$. As the  exact co-distribution $\dd H$ in Definition \ref{d7.3} is generated by a set of co-vector fields,
it is easy to prove that $L_f(\dd h)=\dd L_f(h)$. This fact
 ensures that $L_f(\dd H)\in {\cal E}_*$.
\end{rem}

Recall (\ref{1.4}). Denote by
$H=(h_1(x(t)),\cdots,h_p(x(t)))^{\mathrm{T}}$.
Then it is clear that
\begin{align}\label{7.9}
{\cal H}_*=\Span(\Row(\dd H))
\end{align}
is exact.

By definition, we have the following result.

\begin{prp}\label{p7.4} If the exact co-distribution ${\cal H}_*$ defined by (\ref{7.9}) is $f$ and $g_i$, $i\in [1,m]$ invariant, then there exists an exact OR-system of (\ref{1.4}).
\end{prp}

\noindent{\it Proof.} Assume ${\cal H}_*$ is $f$ and $g_i$, $i\in [1,m]$ invariant, then there exists $\Xi_i\in {\cal M}_{p\times p}$, $i\in[0,p]$, such that
\begin{align*}
\begin{array}{l}
L_f(\dd H)=\Xi_0 \dd H,~L_{g_i}(\dd H)=\Xi_i \dd H,~i\in [1,m].
\end{array}
\end{align*}

Note that $\Xi_i$ are smooth matrices. By Assumption  A3, we have
\begin{align*}
\begin{array}{l}
L_f(H)=\Xi_0 H, \quad L_{g_i}(H)=\Xi_i H,\quad i\in [1,m].
\end{array}
\end{align*}

Then the exact OR-system is
\begin{align*}
\dot{y}=\Xi_0y+\dsum_{i=1}^m\Xi_iu_i.
\end{align*}
\hfill $\Box$

Next, we consider the extended OR-system.
Assume ${\cal H}_*$ is not $f$ and $g_i$, $i\in [1,m]$ invariant. Then we need the extended OR-system. To this end, we calculate the smallest exact co-distribution containing ${\cal H}_*$, called its $\{f,g\}$-closure, and denoted by $\overline{\cal H}^{\{f,g\}}_*$. We need the following algorithm:
\begin{align}\label{7.401}
\begin{array}{l}
{\cal H}^0_*:=\Span\{\dd h_j\;|\;j\in [1,p]\},\\
{\cal H}^{k+1}_*:={\cal H}^k_*+L_f\left({\cal H}^k_*\right)+\dsum_{i=1}^mL_{g_i}\left({\cal H}^k_*\right),\quad k\geq 0.\\
\end{array}
\end{align}
Then obviously when ${\cal H}^{k^*+1}_*={\cal H}^{k^*}_*$, we have 
$\overline{\cal H}^{\{f,g\}}_*={\cal H}^{k^*}_*$.

From the above construction, we have the following result.

\begin{prp}\label{p7.5} Using ${\cal H}^{k^*}_*$ to replace $\cal{H}_*$ yields an exact OR-system, which  is the (smallest) extended OR-system of the original system (\ref{1.4}).
\end{prp}

Finally, we consider the feedback OR-system of (\ref{1.4}).

\begin{dfn}\label{d7.501} A distribution $\D\subset V^*(W)$ is $(f,g)$-invariant, if $\ad_f(\D)\subset \D+{\cal G}$,
where ${\cal G}=\Span(g_1,\cdots,g_m)$.
\end{dfn}

To verify $(f,g)$ invariance, we need the following Quike lemma \cite{isi95}.

\begin{lem}\label{l7.6}
Assume $\D\subset V^*(W)$ is $(f,g)$-invariant, then there exists a neighborhood $0\in U\subset W$ and a smooth mapping
$\a_i(x)\in C^{\infty}(U)$, $i\in [1,m]$, such that $\D$ is $f+\dsum_{i=1}^mg_i\a_i$ invariant.
\end{lem}

Because of Lemma \ref{l7.6}, hereafter we only nee to consider the local case, unless $\D$ is globally $f+\dsum_{i=1}^mg_i\a_i$ invariant.

\begin{lem}\label{l7.7} If $\D^m\subset V^*(W)$ is the largest $(f,g)$-invariant distribution contained in $\ker(H)$, then $\D^m$ is involutive.
\end{lem}

\noindent{\it Proof.} Assume $\D_{m}$ is the largest $(f,g)$-invariant distribution contained in $\ker(H)$, and $X,Y\in \D_{m}$. Then
$$
[X,Y](h)=L_XL_Y(h)-L_YL_X(h)=0,\quad \forall h\in \Row(H).
$$
That is, $[X,Y]\in \D_{m}$.
\hfill $\Box$

\begin{lem}\label{l7.8} Let $\D$ be an $f$-invariant distribution. Denote by
\begin{align*}
\D^{\perp}:=\Span\left\{\dd \xi \;|\; \xi \in C^{\infty}(U), L_X(\xi)=0,\;\forall X\in \D\right\}.
\end{align*}
Then
$\D^{\perp}$ is $f$-invariant.
\end{lem}

\noindent{\it Proof.} Let $\xi\in \D^{\perp}$ and $v\in \D$. Since $\D$ is $f$-invariant, $[f,v]:=v'\in \D$.
Hence $[f,v](\xi)=0$, $\forall \xi\in \D^{\perp}.
$
Now
$$
[f,v](\xi)=L_fL_v(\xi)-L_vL_f(\xi)=L_v(L_f(\xi))=0, \quad  \forall \xi\in \D^{\perp}.
$$
That is $L_f(\xi)\in \D^{\perp}$, $\forall \xi\in \D^{\perp}$,
which means that $\D^{\perp}$ is $f$-invariant.

Using Lemma \ref{l7.8} and a similar argument as in the linear case, we have the following result.

\begin{thm}\label{t7.9}  Let $\D_m$ be the largest $(f,g)$-invariant distribution contained in $\ker({\cal H}_*)$. Then
$\D_m^{\perp}$ is the smallest $(f,g)$-invariant exact co-distribution containing ${\cal H}_*$.
\end{thm}

The algorithm for the largest $(f,g)$-invariant distribution contained in $\ker(H)$ is well known \cite{isi95}. Therefore, $\D_m^{\perp}$ is also easily computable. The feedback OR-system can be constructed as follows:

\begin{prp}\label{p7.10} Using $\D_m^{\perp}$ to replace $H$, the exact OR-system obtained is the (smallest) feedback (extended) OR-system
of the original system (\ref{1.4}).
\end{prp}

Finally, we give a simple example to describe this.

\begin{exa}\label{e7.11}
Consider a nonlinear system
\begin{align}\label{7.10}
\begin{array}{l}
\begin{cases}
\dot{x}_1(t)=x_2(t)+x_3(t)+x_2^2(t),\\
\dot{x}_2(t)=x_1(t)-x_3(t)+x_3^2(t)+u(t),\\
\dot{x}_3(t)=x_3(t)+x_2^2(t),\\
\end{cases}\\
y(t)=x_1(t)-x_3(t).
\end{array}
\end{align}

Using algorithm \ref{7.401}, we have
$$
\begin{array}{l}
G^{\perp}=\Span\left(\begin{bmatrix}1&0&0\end{bmatrix},
\begin{bmatrix}0&0&1\end{bmatrix}
\right);\\
\Omega_0=\Span\left(\begin{bmatrix}1&0&-1\end{bmatrix}
\right);\\
\Omega_k=\Span\left(\begin{bmatrix}1&0&-1\end{bmatrix},
\begin{bmatrix}0&1&0\end{bmatrix}
\right),\quad k\geq 1.\\
\end{array}
$$
Then an efficient feedback can be obtained as
$$
u(t)=-x_3^2+v(t),
$$
The closed-loop system becomes
\begin{align}\label{7.11}
\begin{array}{l}
\dot{x}(t)=\tilde{f}(t)+gv,\\
y(t)=x_1(t)-x_3(t).
\end{array}
\end{align}
where
$$
\tilde{f}(t)=(x_2(t)+x_3(t)+x_2^2(t),x_1(t)-x_3(t),x_2(t)+x_2^2(t))^T.
$$
It follows that  the smallest  $(f,g)$-invariant co-distribution containing observers is
$$
\D_m=\Span\left\{\begin{bmatrix}1&0&-1\end{bmatrix},
\begin{bmatrix}0&1&0\end{bmatrix}
\right\}.
$$
Hence, $y=x_1-x_2$ is not $(f,g)$-invariant. We need to find the $(f,g)$-invariant OR-system.
Consider
$$
\begin{array}{l}
L_{\tilde{f}}h(x)=x_2,\quad
L_g(h(x))=0.
\end{array}
$$
we know that
$$
\D_m=\Span\{\dd y,\dd x_2\}.
$$
Hence the extended $(f,g)$-invariant OR containing $y$ is
\begin{align}\label{7.12}
\begin{array}{l}
\dot{y}(t)=x_2(t),\\
\dot{x}_2(t)=y(t)+v(t).
\end{array}
\end{align}
\end{exa}

\begin{rem}\label{r7.12}
In Example \ref{e7.11}, note that  \eqref{7.10} is not a state-feedback linearizable system, while the extended feedback OR-system is   feedback linearizable. On the other hand, if the original system is linearizable, its extended OR-system is also linearizable. This suggests that the OR-system approach may offer advantages in the analysis of nonlinear systems.
\end{rem}

\section{Conclusion}

This paper presents a technique for constructing the OR-systems. First, by employing cross-dimensional projection and DK-STP, two types of approximate OR-systems are proposed. Subsequently, we address the conditions under which an OR-system is exact, meaning that it exhibits precisely the same dynamics as the observers of the original system. The methodology for constructing an exact OR-system is also provided.

When it is not possible to construct an exact OR-system, we propose algorithms for constructing an extended OR-system, in which the observers of the original system are included as part of the state variables. Furthermore, the (minimal) state feedback (extended) OR-system is developed. Finally, all the aforementioned results for linear control systems are extended to affine nonlinear control systems. The approaches for extended and feedback OR-systems in continuous-time control systems are also applicable to discrete-time control systems.

This work aims to propose a method for model reduction of complex systems. When dealing with large-scale complex systems, the entire state space may be too complicated to analyze or even to model. In such cases, we can focus on certain observed phenomena, represented by the OR-system of the original system. By studying various OR-systems derived from the original system, we are able to gain insight into, and potentially manipulate, the underlying complex system.

Essentially, OR-realization can be regarded as a novel model reduction technique. There remain many related problems for further investigation, some of which are listed below:

 \begin{itemize}
 \item[(i)] For approximate OR-realization, how can the approximation error be estimated? Furthermore, how can the error be reduced to meet application requirements?
 \item[(ii)] What constitutes the optimal OR-approximation for a given system? In particular, the minimal dimension of the best approximation remains an open question.
 \item[(iii)] The OR-systems technique has significant potential in model reduction in large-scale complex systems. In particular, its application to dimension-varying systems \cite{che23} brings notable challenges.
 \end{itemize}

\section*{Appendix:~Calculation of $(A,B)$-Invariant Subspaces (for Example \ref{e6.2.11})}

First, we recall the following algorithm (see \cite{won79}, page 95):
\begin{alg}\label{aa.2.1}

\begin{itemize}
\item[(i)]
Set ${\cal V}_0=X$.
\item[(ii)] For $k\geq 1$, compute
${\cal V}_k=X\cap {\cal A}^{-1}\left({\cal B}+{\cal V}_{k-1}\right).
$
\item[(iii)] If ${\cal V}_{k^*+1}={\cal V}_{k^*}$, stop.
\end{itemize}
The largest $(A,B)$-invariant subspace contained in $X$ is then given by ${\cal V}_{k^*}$.
\end{alg}

Note that in the above algorithm
${\cal A}^{-1}({\cal S}):=\{x\in \R^n\;|\;Ax\in {\cal S}\}.
$

To compute ${\cal A}^{-1}({\cal S})$, we use the following lemma.

\begin{lem}\label{la.2.2} ${\cal A}^{-1}({\cal S})=\left[A^\mathrm{T} {\cal S}^{\perp}\right]^{\perp}$, where ${\cal A}=\Span(\Col(A))$.
\end{lem}

\noindent{\it Proof.}
Let $x\in A^\mathrm{T}{\cal S}^{\perp}$. Then there exists $y\in {\cal S}^{\perp}$, which means $s^\mathrm{T}y=0$, $\forall s\in {\cal S}$, such that $x=A^\mathrm{T}y$. Let $z\in {\cal A}^{-1}({\cal S})$, i.e., $Az=s'\in {\cal S}$. So, $z^\mathrm{T}x=z^\mathrm{T}A^\mathrm{T}y=(s')^\mathrm{T}y=0$. That is, $x\in ({\cal A}^{-1}{\cal S})^{\perp}$. Hence,
$\left(A^\mathrm{T}{\cal S}^{\perp}\right)^{\perp}\supset {\cal A}^{-1}({\cal S})$.

Conversely, suppose $x\in (A^\mathrm{T}{\cal S}^{\perp})^{\perp}$. Then $x^\mathrm{T}A^\mathrm{T}y=0$, $\forall y\in {\cal S}^{\perp}$, which implies
$Ax\in ({\cal S}^{\perp})^{\perp}={\cal S}$. Therefore, $x\in {\cal A}^{-1}({\cal S})$, and thus $(A^\mathrm{T}{\cal S}^{\perp})^{\perp}\subset {\cal A}^{-1}({\cal S}).
$
\hfill $\Box$

Applying Algorithm \ref{aa.2.1}, we set
${\cal V}_0:={\cal C}^{\perp}=\Span(\Col(V_0))$,
where $
V_0=\left[\begin{smallmatrix}
-1&0&0&0&0\\
0&1&0&0&0\\
0&0&1&0&0\\
0&0&0&1&0\\
1&0&0&0&0\\
0&0&0&0&1\\
\end{smallmatrix}\right],~{\cal S}_0={\cal B}+{\cal V}_0={\cal V}_0.
$
\begin{align*}
    {\cal S}^{\perp}_0=\Span\left(\left[\begin{smallmatrix}1&0&0&0&1&0\end{smallmatrix}\right]^{\mathrm{T}}\right),\\
A^\mathrm{T}{\cal S}^{\perp}_0=\Span\left(\left[\begin{smallmatrix}0&1&0&0&1&0\end{smallmatrix}\right]^{\mathrm{T}}\right).
\end{align*}
\begin{align*}
    {\cal A}^{-1}({\cal S}_0)&=\left(A^\mathrm{T}{\cal S}_0^{\perp}\right)^{\perp}=\Span\left(\Col\left(\left[\begin{smallmatrix}
1&0&0&0&0\\
0&0&0&1&0\\
0&1&0&0&0\\
0&0&1&0&0\\
0&0&0&-1&0\\
0&0&0&0&1\end{smallmatrix}\right]\right) \right).
\end{align*}
\begin{align*}
    {\cal V}_1={\cal V}_0\bigcap {\cal A}^{-1}({\cal S}_0)=\Span\left(\Col\left(\left[\begin{smallmatrix}
1&0&0&0\\
1&0&0&0\\
0&1&0&0\\
0&0&1&0\\
-1&0&0&0\\
0&0&0&1\end{smallmatrix}\right]\right) \right).
\end{align*}
$$
{\cal S}_1={\cal B}+{\cal V}_1={\cal V}_1.
$$
$$
{\cal S}^{\perp}_1=\Span\left(\Col\left(\left[\begin{smallmatrix}
1&1\\
-1&1\\
0&0\\
0&0\\
0&2\\
0&0\end{smallmatrix}\right]\right) \right).
$$
$$
A^\mathrm{T}{\cal S}^{\perp}_1=\Span\left(\Col\left(\left[\begin{smallmatrix}
0&0\\
1&1\\
-1&1\\
0&0\\
0&2\\
0&0\end{smallmatrix}\right]\right) \right).
$$
$$
{\cal A}^{-1}({\cal S}_1)=\left(A^\mathrm{T}{\cal S}_1^{\perp}\right)^{\perp}
=\Span\left(\Col\left(\left[\begin{smallmatrix}
1&0&0&0\\
0&0&0&1\\
0&0&0&1\\
0&1&0&0\\
0&0&0&-1\\
0&0&1&0\end{smallmatrix}\right]\right) \right).
$$
$$
{\cal V}_2={\cal V}_1\bigcap {\cal A}^{-1}({\cal S}_1)
=\Span\left(\Col\left(\left[\begin{smallmatrix}
0&0&1\\
0&0&1\\
0&0&1\\
1&0&0\\
0&0&-1\\
0&1&0\end{smallmatrix}\right]\right) \right).
$$
$$
{\cal S}_2={\cal B}+{\cal V}_2.
$$
It is straightforward to verify that
${\cal S}_2={\cal S}_1$.
Therefore,
${\cal V}_3={\cal V}_2$,
which is the largest $(A,B)$-invariant subspace contained in ${\cal C}_*^{\perp}$.

\end{document}